\documentclass[12pt]{amsart}

\usepackage{graphicx}
\usepackage[leqno]{amsmath}
\usepackage{amssymb}
\usepackage[all]{xy}
\usepackage{graphicx}
\usepackage[usenames,dvipsnames]{color}
\usepackage[normalem]{ulem}
\usepackage{tikz-cd}

\newtheorem{thm}{Theorem}[section]
\newtheorem{notation}[thm]{Notation}
\newtheorem{lemma}[thm]{Lemma}
\newtheorem{cor}[thm]{Corollary}
\newtheorem{prop}[thm]{Proposition}

\newtheorem{remark}[thm]{Remark}

\theoremstyle{definition}
\newtheorem{definition}[thm]{Definition}

\date{\today}

\begin{document}

\title[$C^0$ approximations of foliations]
{$C^0$ approximations of foliations}

\title[$C^0$ approximations of foliations]
{$C^0$ approximations of foliations}

\author[Kazez]{William H.  Kazez}
\address{Department of Mathematics, University of Georgia, Athens, GA 30602}
\email{will@math.uga.edu}

\author[Roberts]{Rachel Roberts}
\address{Department of Mathematics, Washington University, St.  Louis, MO 63130}
\email{roberts@math.wustl.edu}

\thanks{This work was partially supported by grants from the Simons Foundation (\#244855 to William Kazez, \#317884 to Rachel Roberts)}



\begin{abstract} Suppose that $\mathcal F$ is a transversely oriented, codimension one foliation of a connected, closed, oriented 3-manifold.  Suppose also that $\mathcal F$ has continuous tangent plane field and is {\sl taut}; that is, closed smooth transversals to $\mathcal F$ pass through every point of $M$.  We show that if $\mathcal F$ is not the product foliation $S^1\times S^2$, then $\mathcal F$ can be $C^0$ approximated by weakly symplectically fillable, universally tight, contact structures.  This extends work of Eliashberg-Thurston on approximations of taut, transversely oriented $C^2$ foliations to the class of foliations that often arise in branched surface constructions of foliations.  This allows applications of contact topology and Floer theory beyond the category of $C^2$ foliated spaces.
\end{abstract}

\maketitle

\section{Introduction}

In \cite{ET}, Eliashberg and Thurston introduce the notion of {\it confoliation} and prove that when $k\ge 2$, a transversely oriented $C^k$~foliation $\mathcal F$ of a closed, oriented 3-manifold not equal to $S^1\times S^2$ can be $C^0$~approximated by a pair of $C^{\infty}$ contact structures, one positive and one negative.  They also prove that when $\mathcal F$ is also taut, any contact structure sufficiently close to the tangent plane field of $\mathcal F$ is weakly symplectically fillable and universally tight.  

The focus of this paper is $C^{k,0}$ foliations (see the next section for definitions).  These foliations are less smooth than those studied by Eliashberg and Thurston.  In the context of applications of $C^{k,0}$ foliations, the natural definition of taut is a foliation for which closed smooth transversals pass through every point of the manifold.

\begin{definition} A foliation $\mathcal F$ of $M$ is ${\sl taut}$ if for each $p\in M$, there is a simple closed curve everywhere transverse to $\mathcal F$ and passing through $p$.
\end{definition}

Given sufficient smoothness, this is equivalent to the usual definition which stipulates the existence of transversals through every leaf of the foliation \cite{KR5}.  In addition when dealing with $C^{k,0}$ foliations, it is important to specify whether smooth or topological transversals are being used \cite{KR5}.  We use the terms {\sl transverse, transversal,} and {\sl transversely} in the smooth sense; that is, they refer to smooth objects intersecting so that the associated tangent spaces intersect minimally.  In contrast, a curve is {\sl topologically transverse} to $\mathcal F$ if no nondegenerate subarc   is isotopic, relative to its endpoints, into a leaf of $\mathcal F$.

In this paper, we complete the project started in \cite{KR2} and prove that the requirement that $\mathcal F$ be $C^2$ can be weakened to the condition that $\mathcal F$ be $C^{1,0}$; equivalently, that the tangent plane field of $\mathcal F$ be defined and continuous.  Our main theorem is the following.

\begin{thm}\label{maincor}
Any taut, transversely oriented $C^{1,0}$ foliation on a closed, connected, oriented 3-manifold $M\ne S^1\times S^2$ can be $C^0$ approximated by both a positive and a negative smooth contact structure, $\xi^+$ and $\xi^-$.  These contact structures $(M,\xi^+)$ and $(-M,\xi^-)$ are weakly symplectically fillable and universally tight.
\end{thm}

In addition, we take the proof of this theorem as an opportunity to revisit classical foliation results in the topological setting, and attempt to state as clearly as possibly what is true in the TOP category.  Proofs are given when we could not find them in the literature.

Tautness of a $C^{1,0}$ foliation $\mathcal F$ guarantees the existence of a transverse, smooth, volume preserving flow (Theorem~6.1 of \cite{KR5}).  The existence of a volume preserving flow transverse to $\mathcal F$ is used in Proposition~3.2.2 of \cite{ET} to conclude weak symplectic fillability of the approximating contact structure.  Tightness then follows from Theorem~3.2.4, and universal tightness by Proposition~3.5.5, of \cite{ET}.  Theorem~\ref{maincor} therefore implies the following.  

\begin{cor} 
Let $\mathcal F$ be a transversely oriented $C^{1,0}$ foliation on a closed, connected, oriented 3-manifold $M\ne S^1\times S^2$.  If $\Phi$ is a smooth volume preserving flow transverse to $\mathcal F$, then any positive contact structure transverse to $\Phi$ is weakly symplectically fillable and universally tight.
\end{cor}

Jonathan Bowden \cite{bowden} has also independently proved Theorem~\ref{maincor}, and hence any accompanying corollaries, when $\mathcal F$ is $C^{\infty,0}$ and not a foliation of $T^3$ by planes.  His approach is similar to ours but uses different propagation techniques.  Bowden also makes the interesting observation that the condition that the foliation be taut is necessary only to conclude that there is a volume preserving transverse flow, and hence that the approximating contact structures are weakly symplectically fillable and universally tight.  There always exists a smooth transverse flow, and this is all that is necessary to find the approximating contact structures.  This is true of our proof also.  In summary, we have: 

\begin{thm}[Theorem~1.2, \cite{bowden}] Let $M$ be a closed, connected, oriented 3-manifold, and let $\mathcal F$ be a transversely oriented $C^{1,0}$ foliation on $M$.  Then $\mathcal F$ can be $C^0$ approximated by a positive (respectively, negative) contact structure if and only if $\mathcal F$ is not a foliation of $S^1\times S^2$ by spheres.  
\end{thm}

In this paper $M$ will denote a connected, closed, oriented 3-manifold, and $\mathcal F$ will be a transversely oriented (not necessarily taut) foliation in $M$.  

There are two useful notions that we use in describing approximations of a given foliation $\mathcal F$ by either another foliation, a contact structure, or indeed any $2$-plane field $\xi$ on $M$.  The first is $C^0$ approximation.  This can be defined using the standard topology, or a compatible metric, on the Grassmannian bundle of continuous tangent 2-planes over $M$.  Pick a compatible metric.  If for all $\epsilon > 0$ there exists a 2-plane field $\xi$ on $M$ of a particular type (e.g., corresponding to a foliation in Theorem~\ref{main}, or a contact structure in Corollary~\ref{maincor}) that is within $\epsilon$ of $T{\mathcal F}$, then we say $\mathcal F$ is $C^0$ close to a 2-plane bundle of the type of $\xi$.

The starting point for the second notion of approximation is a foliation $\mathcal F$ together with a transverse flow $\Phi$.  A tangent 2-plane distribution $\xi$ on $M$ is {\sl $\Phi$-close} to $\mathcal F$ if it is also transverse to $\Phi$.  Since the notion of $\Phi$-closeness is purely topological, it is very well suited to the study of continuous plane fields.  When $\mathcal F$ is taut, there exists a volume preserving flow $\Phi$ transverse to $\mathcal F$, and this notion of $\Phi$ approximation is particularly useful, and is the focus of \cite{KR2}, because it is sufficient that a contact structure $\xi$ be $\Phi$-close to $\mathcal F$ in order to conclude that it is weakly symplectically fillable and universally tight.  Clearly $\Phi$-close is implied by $C^0$ close, and it is the latter notion that is the focus of this paper.

The construction in \cite{ET} of a contact structure approximating a foliation consists of two steps.  First, a contact structure is constructed in neighborhoods of curves in leaves of the foliation about which there is contracting holonomy.  Next the foliation is used to propagate the contact structure to the remainder of the manifold.

To carry out this strategy for $C^{1,0}$ foliations, extra care is required.  The preferred neighborhoods of curves must be chosen to be particularly thin, and it may be that a foliation has no such thin neighborhoods.  In such a case, a new foliation that $C^0$ approximates the first and has thin, contracting holonomy is constructed.  This is done by the method of generalized Denjoy blow up and is described in Section~\ref{Denjoy}.

The key issue in attempting to propagate a contact structure defined only on a subset $V$ of the manifold to the entire manifold is the possibility that not every point outside of $V$ is connected to $V$ by a leaf of the foliation.  The main ideas in finding such a set $V$ are discussed in Section~\ref{minimal}.  Since we are changing the foliation and studying a possibly changing minimal set, arguments have to be made that this is a finite procedure, and that the desired $V$ is just a finite union of thin contracting neighborhoods.

The culmination of this strategy is Theorem~\ref{metriccreateholonbhds}, which can be roughly stated as

\begin{thm}\label{main} A transversely oriented $C^{1,0}$ foliation $\mathcal F$ not equal to $S^1\times S^2$ is $C^0$ close to a $C^{\infty,0}$ foliation $\mathcal F'$ for which there exists a subset $V$ of $M$ with the following properties.  First, $V$ is a disjoint union of finitely many thin solid tori each of which is a standard neighborhood of a closed curve in a leaf for which the holonomy has a contracting interval.  Second, $M$ is $V$-transitive with respect to $\mathcal F'$.
\end{thm}

Much of the focus of the paper is on creating and approximating nontrivial holonomy in foliations.  Section~\ref{noholo} covers the case in which no nontrivial holonomy exists, including the case that all leaves of the foliation are planes.  Throughout the paper, {\sl plane} is used in the topological sense; that is, a plane is a surface homeomorphic to $\mathbb R^2$.  We use $I$ to mean nondegenerate closed interval in places, and to mean $I=[0,1]$ in other places, and the meaning should be clear from context.  

In Section~\ref{propagate} we briefly recall the techniques introduced in \cite{KR2} to propagate a contact structure along a $C^{\infty,0}$ foliation and thereby conclude our main result, Theorem~\ref{maincor}, from Theorem~\ref{main}.

In \cite{calegari}, Calegari proves the theorem, proposed in \cite{Ga93} as folklore in need of proof, that any $C^0$ foliation can be isotoped to a $C^{\infty,0}$ foliation.  This leads to an existence, as opposed to an approximation, result.  That is, the existence of a taut, oriented $C^0$ foliation is sufficient to guarantee the existence of a pair of contact structures, one positive and one negative.  

\begin{cor} Suppose $M$ is a closed oriented 3-manifold that contains a taut oriented $C^0$ foliation $\mathcal F$ not equal to $S^1\times S^2$.  Then $M$ contains a pair of contact structures $\xi^+$ and $\xi^-$, one positive and one negative, that may be chosen arbitrarily $C^0$ close to each other.  These contact structures $(M,\xi^+)$ and $(-M,\xi^-)$ are weakly symplectically fillable and universally tight.
\end{cor}

In \cite{OS}, Ozsvath-Szabo use the Eliashberg-Thurston existence theorem to show that $L$-spaces do not admit transversely orientable, taut, $C^2$ foliations.  Since many constructions of foliations, \cite{DL, DR, G1, g1, g2, g3, G2, G3, Ga1, KRo, KR2, Li, Li3, LR, R, R1, R2}, are not $C^2$, it is useful to be able to remove the smoothness assumption from their theorem.

A $C^0$ foliation is {\sl topologically taut} if there is a topological transversal through every leaf of the foliation.  By Corollary~5.6 of \cite{KR5}, a topologically taut $C^0$ foliation is isotopic to a taut $C^{\infty,0}$ foliation, and hence the differences between the versions of tautness are unimportant when working with foliations up to topological conjugacy.

\begin{cor} An $L$-space does not admit a transversely orientable, topologically taut, $C^0$ foliation.
\end{cor}

One of the motivations for our work is the uniqueness theorem proved by Vogel \cite{vogel}.  He shows that with mild assumptions on the leaves of a $C^2$ foliation, that are necessarily satisfied in an atoroidal, irreducible manifold, sufficiently close approximating contact structures must be isotopic to each other.  Although our work does not address the uniqueness question for approximations of $C^{1,0}$ foliations, one of our main tools, flow box decompositions, seems well suited for comparing a pair of close contact structures.

This suggests questions related to uniqueness that can be asked in an atoroidal, irreducible manifold.
\begin{enumerate}

\item If $\xi_1$ and $\xi_2$ are sufficiently close to a taut, transversely oriented $C^{1,0}$ foliation $\mathcal F$, are they isotopic?

\item Can the contact invariant of an approximating contact structure be computed, and shown not to vanish, directly from the foliation?

\item If $\xi_1$ and $\xi_2$ are sufficiently close to a taut, transversely oriented $C^{1,0}$ foliation $\mathcal F$, are the contact invariants determined by $\xi_1$ and $\xi_2$ necessarily equal?
\end{enumerate}

We thank Vincent Colin, Larry Conlon and John Etnyre for many helpful conversations.  We also thank Jonathan Bowden for helpful feedback on a preliminary version of this paper.

\section{Background}\label{background}

We begin by defining foliations in 3-manifolds with empty boundary.  Near the end of this section, we extend these definitions to 3-manifolds with nonempty boundary that are smooth or smooth with corners; namely, manifolds locally modeled by open sets in $[0,\infty)^3$.

\begin{definition}\label{folndefn1} Let $M$ be a smooth 3-manifold with empty boundary.  Let $k$ and $l$ be non-negative integers or infinity with $l \leq k$.  Both $C^k$ and $C^{k,l}$ {\sl codimension one foliations} $\mathcal F$ are decompositions of $M$ into a disjoint union of $C^k$ immersed connected surfaces, called the {\sl leaves} of $\mathcal F$, together with a collection of charts $U_i$ covering $M$, with $\phi_i:\mathbb R^2 \times \mathbb R \to U_i$ a homeomorphism, such that the preimage of each component of a leaf intersected with $U_i$ is a horizontal plane.  

The foliation $\mathcal F$ is $C^k$ if the charts $(U_i,\phi_i)$ can be chosen so that each $\phi_i$ is a $C^k$ diffeomorphism.  

The foliation $\mathcal F$ is $C^{k,l}$ if for all $i$ and $j$,

\begin{enumerate}

\item the derivatives $\partial_x^{\hskip .015 truein a}\partial_y^{\hskip .015 truein b}\partial_z^{\hskip .015 truein c}$, taken in any order, on the domain of each $\phi_i$ and each transition function $\phi_j^{-1}\phi_i$ are continuous for all $a + b \leq k$, and $c \leq l$, and

\item if $l \geq1$, $\phi_i$ is a $C^1$ diffeomorphism.

\end{enumerate}
\end{definition}

\begin{remark} The smoothness conditions on both the charts and the transition functions are to ensure that the smooth structure on the leaves is compatible with the smooth structure on $M$.  
\end{remark}

In particular, $T\mathcal F$ exists and is continuous if and only if $\mathcal F$ is $C^{1,0}$.  Also notice that $C^{k,l}$ foliations are $C^l$, but not conversely.

Two $C^{k,0}$ foliations $\mathcal F$ and $\mathcal G$ of $M$ are called {\sl $C^{k,0}$ equivalent} if there is a self-homeomorphism of $ M$ that maps the leaves of $\mathcal F$ to the leaves $\mathcal G$, and is $C^k$ when restricted to any leaf of $\mathcal F$.

\begin{remark} Definition~\ref{folndefn1} is an amalgamation of two definitions due to Candel and Conlon (1.2.22, 1.2.24) in \cite{CC}.  In one of their definitions, they allow for the possibility that the ambient manifold is not given a differentiable structure or that it may have a differentiable structure that does not contain $T\mathcal F$ as a subbundle.  Since a topological 3-manifold admits a unique smoothness structure \cite{Mo}, we forego this generality and require leaves of $\mathcal F$ to be $C^k$ immersed in $M$.
\end{remark}

Given a codimension one foliation $\mathcal F$, it is useful to fix a flow $\Phi$ transverse to $\mathcal F$.  Even when the leaves of the foliation are only topologically immersed, Siebenmann shows (Theorem~6.26, \cite{Siebenmann}; see also Chapter IV of \cite{HH}) that there is a 1-dimensional transverse foliation.  For foliations with smoother leaves it is much easier to construct a transverse flow.  We next recall some basic facts about flows.

\begin{definition} [See, for example, Section 17 of \cite{Lee}.] A {\sl global flow} on $M$ is an action of $\mathbb R$ on $M$; that is, a continuous map
$$\phi: M \times \mathbb R \to M$$
such that 
\begin{enumerate}
\item $\phi(\phi(p,s),t) = \phi(p,t+s)$ and
\item $\phi(p, 0)=p$.
\end{enumerate}
For each $p\in M$, there is a curve $\phi_p:\mathbb R\to M$ given by $\phi_p(t) =\phi(p, t)$.  
\end{definition}

\begin{prop} [Proposition 17.3 of \cite{Lee}]
Let $\phi: M \times \mathbb R \to M$ be a smooth global flow.  The vector field $V:M\to TM$ given by $V(p)=\phi_p'(0)$ is smooth, and each curve $\phi_p$ is an integral curve of $V$.  \qed
\end{prop}

We call $V$ the vector field {\sl determined by} the flow $\phi$.  Let $\Phi$ denote the 1-dimensional smooth foliation with leaves the integral curves $\phi_p$ of $V$.  Conversely, since $M$ is compact, any smooth vector field determines a smooth flow.

\begin{thm}[Theorem 17.8 and 17.11 of \cite{Lee}]
Given a smooth vector field $V$ on $M$, there is a unique global flow $\phi: M \times \mathbb R \to M$ such that $V$ is the vector field determined by $\phi$.  \qed
\end{thm}

So any choice of nowhere zero tangent vector field to a smooth 1-dimensional foliation determines a smooth global flow on $M$.

If there is a continuously varying vector field transverse to the leaves of a foliation $\mathcal F$, then $\mathcal F$ is {\sl transversely orientable}.

\begin{thm} [Theorem~6.26, \cite{Siebenmann}, Theorems 1.1.2 and 1.3.2, \cite{HH}] \label{perp} Let $\mathcal F$ be a codimension one, transversely oriented $C^0$ foliation.  There is a continuous flow $\Phi$ transverse to $\mathcal F$.  

When $M$ is $C^{1,0}$, this flow can be chosen to be smooth.  In fact, there is a smooth nowhere zero vector field $V:M\to TM$ everywhere positively transverse to $\mathcal F$, and if $M$ is given a Riemannian metric, then $V$ can be chosen to lie arbitrarily $C^0$ close to the continuous vector field of vectors perpendicular to $T\mathcal F$.
\end{thm}

\begin{proof} We are primarily interested in the case that $\mathcal F$ is $C^{1,0}$.  Since in this case the proof is both immediate and enlightening, we reproduce it here.  Fix a Riemannian metric on $M$.  Since $\mathcal F$ is transversely oriented, each tangent plane $T_p\mathcal F$ has a preferred orientation.  For each $p\in M$, let $V_\perp(p)$ denote the positive unit normal to the tangent plane $T_p\mathcal F$.  Approximate $V_\perp(p)$ by a smooth vector field $V$.  If the approximation is taken close enough it will be non-zero and transverse to $T\mathcal F$.
\end{proof}

When a foliation $\mathcal F$ and a transverse flow $\Phi$ are understood, a submanifold of positive codimension in $M$ is called {\sl horizontal} if each component is a submanifold of a leaf of $\mathcal F$ and {\sl vertical} if it can be expressed as a union of subsegments of the flow $\Phi$.  A codimension-0 submanifold $X$ of $M$ is called {\sl $(\mathcal F,\Phi)$-compatible} if its boundary is piecewise horizontal and vertical, and hence $\mathcal F$ and $\Phi$ restrict naturally to foliation and flow on $X$.  If $X$ is $(\mathcal F,\Phi)$-compatible, let $\partial_v X$ denote its vertical boundary and let $\partial_h V$ denote its horizontal boundary.

\begin{definition} \cite{KR2} \label{flowboxdefn} 
Let $\mathcal F$ be either a $C^k$ or $C^{k,l}$ foliation, and let $\Phi$ be a transverse flow.  A {\sl flow box}, $F$, is an $(\mathcal F,\Phi)$ compatible closed chart, possibly with corners.  That is, it is a submanifold diffeomorphic to $D\times I$, where $D$ is either a closed $C^k$ disk or polygon (a closed disk with at least three corners), $\Phi$ intersects $F$ in the arcs $\{(x,y)\}\times I$, and each component of $D\times \partial I$ is embedded in a leaf of $\mathcal F$.  The components of $\mathcal F\cap F$ give a family of graphs over $D$.  

In the case that $D$ is a polygon, it is often useful to view the disk $D$ as a 2-cell with $\partial D$ the cell complex obtained by letting the vertices correspond exactly to the corners of $D$.  Similarly, it is useful to view the flow box $F$ as a 3-cell possessing the product cell complex structure of $D\times I$.  Then $\partial_h F$ is a union of two (horizontal) 2-cells and $\partial_v F$ is a union of $c$ (vertical) 2-cells, where $c$ is the number of corners of $D$.  In the case that $D$ has no corners, we abuse language slightly and consider $\partial_h F$ to be a union of two (horizontal) 2-cells and $\partial_v F$ to be a single vertical face, where the face is the entire vertical annulus $\partial D\times I$.

Suppose $V$ is either empty or else a compact, $(\mathcal F,\Phi)$ compatible, codimension-0 submanifold of $M$.  A {\sl flow box decomposition of $M$ rel $V$}, or simply {\sl flow box decomposition of $M$}, if $V=\emptyset$, is a decomposition of $M\setminus \text{int} V$ as a finite union $M = V\cup (\cup_iF _i)$ where

\begin{enumerate}
\item each $F _i$ is a flow box,

\item the interiors of $F _i$ and $F _j$ are disjoint if $i \neq j$, and 

\item if $i \ne j$ and $F_i \cap F_j $ is nonempty, it must be homeomorphic to a point, an interval, or a disk that is wholly contained either in $\partial_h F_i \cap \partial_h F_j$ or in a single face in each of $\partial_v F_i$ and $\partial_v F_j$.

\end{enumerate}
\end{definition}

In \cite{KR4}, we develop a theory of flow box decompositions and show that they are particularly well suited to the study of codimension one foliations.  Their role is similar to that played by triangulations and branched surfaces, but they are perhaps better suited to the consideration of differentiability properties.

\begin{lemma} \cite{KR2} \label{fbdexist}
Let $\mathcal F$ be a $C^0$ foliation and $\Phi$ a $C^0$ transverse flow.  There is a flow box decomposition for $(M,\mathcal F,\Phi)$.  When $\mathcal F$ is $C^{k,0}$, $k\ge 1$, and $\Phi$ is smooth, this flow box decomposition can be chosen to be $C^k$.  \qed
\end{lemma}

Two 2-plane bundles, for instance two contact structures, are said to be {\sl $C^0$ close}, if at each point, the associated 2-planes are close in the associated Grassmann bundle of 2-planes.  Two $C^{1,0}$ foliations, $\mathcal F$ and $\mathcal G$, are {\sl $C^0$ close} if the associated two plane bundles $T\mathcal F$ and $T\mathcal G$ are $C^0$ close.  A two plane bundle, for instance a contact structure, is {\sl $C^0$ close} to a $C^{1,0}$ foliation $\mathcal F$ if it is $C^0$ close to $T\mathcal F$.  A diffeomorphism $C^1$ close to the identity map will preserve $C^0$ proximity of foliations and contact structures.

As the next theorem and its corollary show, there is often no loss of generality in restricting attention to foliations with smooth leaves.

\begin{thm}\cite{calegari,KR4}\label{calegarikr}
Suppose $\mathcal F$ is a $C^{1,0}$ foliation in $M$.  Then there is an isotopy of $M$ taking $\mathcal F$ to a $C^{\infty,0}$ foliation $\mathcal G$ which is $C^0$ close to $\mathcal F$.  If $\Phi$ is a smooth flow transverse to $\mathcal F$, the isotopy may be taken to map each flow line of $\Phi$ to itself.  
\end{thm}

\begin{cor} If every $C^{\infty,0}$ taut transversely oriented foliation can be $C^0$ approximated by positive and negative contact structures, then the same is true for every $C^{1,0}$ taut transversely oriented foliation.  \qed
\end{cor}

\section{Holonomy neighborhoods}

Let $\gamma$ be an {\sl oriented} simple closed curve in a leaf $L$ of $\mathcal F$, and let $p$ be a point in $\gamma$.  We are interested in the behavior of $\mathcal F$ in a neighborhood of $\gamma$.  Let $h_{\gamma}$ be a holonomy map for $\mathcal F$ along $\gamma$, and let $\sigma$ and $\tau$ be small closed segments of the flow $\Phi$ which contain $p$ in their interiors and satisfy $h_{\gamma}(\tau) = \sigma $.  Choose $\tau$ small enough so that $\sigma \cup \tau$ is a closed segment and not a loop.  Notice that $\sigma\cap\tau$ is necessarily a closed segment containing $p$ in its interior.  There are three possibilities:
\begin{enumerate}
\item $\sigma=\tau$,
\item one of $\sigma$ and $\tau$ is properly contained in the other, or
\item $\sigma\cap\tau$ is properly contained in each of $\sigma$ and $\tau$.  
\end{enumerate}

We will need to consider very carefully a regular neighborhood of $\gamma$ which lies nicely with respect to both $\mathcal F$ and $\Phi$.  To this end, restrict attention to foliations $\mathcal F$ which are $C^{\infty,0}$ and transverse flows $\Phi$ which are smooth, and suppose that $\gamma$ is smoothly embedded in $L$.  Let $A$ be the closure of a smooth regular neighborhood of $\gamma$ in $L$; so $A$ is a smoothly embedded annulus in $L$.

\begin{lemma}\label{holnbddefn} Suppose $\mathcal F$ is $C^{\infty,0}$ and $\Phi$ is smooth.  If $\tau$ and $A$ are chosen to be small enough, there is a compact submanifold $V$ of $M$, smoothly embedded with corners, which satisfies the following:
\begin{enumerate}
\item $V$ is homeomorphic to a solid torus,
\item $\partial V$ is piecewise vertical and horizontal; namely, $\partial V$ decomposes as a union of subsurfaces $\partial_v V\cup \partial_h V$, where $\partial_v V$ is a union of flow segments of $\Phi$ and $\partial_h V$ is a union of two surfaces $L_-$ and $L_+$, each of which is either a disk or an annulus, contained in leaves of $\mathcal F$,
\item each flow segment of $\Phi\cap V$ runs from $L_-$ to $L_+$,
\item $\tau$ is a component of the flow segments of $\Phi\cap V$, and
\item $A$ is a leaf of the foliation $\mathcal F\cap V$.
\end{enumerate}
\end{lemma}

\begin{proof}
Cover a small open neighborhood of $\gamma$ by finitely many smooth flow boxes.  By passing to a smaller $\tau$ and $A$ as necessary, we may suppose that $A$ is covered by two flow boxes with union, $V$, satisfying the properties (1)-(5).
\end{proof}

\begin{notation}\label{notationattracting}
Denote the neighborhood $V$ of Lemma~\ref{holnbddefn} by $V_{\gamma}(\tau,A)$.  
\end{notation}

Notice that if $\tau=\sigma$, then $V_{\gamma}(\tau,A)$ is diffeomorphic to $A\times I$, where $I$ is a nondegenerate closed interval.  Otherwise, there is a unique smooth vertical rectangle, $R$ say, so that the result of cutting $V_{\gamma}(\tau,A)$ open along $R$, and taking the metric closure, is diffeomorphic to a solid cube.

\begin{notation} \label{notationQ}
Let $R_{\gamma}(\tau,A)$ denote any smooth vertical rectangle embedded in $V_{\gamma}(\tau,A)$ such that the result of 
cutting $V_{\gamma}(\tau,A)$ open along $R$, and taking the metric closure, is diffeomorphic to a solid cube.  When $V_{\gamma}(\tau,A)$ is not diffeomorphic to a product, $R_{\gamma}(\tau,A)$ is uniquely determined.  Let $Q_{\gamma}(\tau,A)$ denote the resulting solid cube; so

$$V_{\gamma}(\tau,A) | {R_{\gamma}(\tau,A)}=Q_{\gamma}(\tau,A).$$
\end{notation}
Note that if $\gamma$ is essential, then $Q_{\gamma}(\tau,A)$ can be viewed as a $(\tilde{\mathcal F},\tilde{\Phi})$ flow box, where $(\tilde{\mathcal F},\tilde{\Phi})$ is the lift of $(\mathcal F,\Phi)$ to the universal cover of $M$.

\begin{definition} The neighborhood $V_{\gamma}(\tau,A)$ is called the {\sl holonomy neighborhood determined by $(\tau,A)$}, and is called an {\sl attracting neighborhood} if $h_\gamma(\tau)$ is contained in the interior of $\sigma$.  
\end{definition}

Notice that at most one of $V_{\gamma}(\tau)$ and $V_{-\gamma}(\tau)$ can be attracting.  More generally, consider the fixed point set of $h_{\gamma}$, $\mbox{Fix}(h_{\gamma})\subset \tau$.  $\mbox{Fix}(h_{\gamma})$ is closed and cuts $\tau$ into open intervals on which $h_{\gamma}$ is either strictly increasing or strictly decreasing.  Identify $(\tau, p)$ with $(I,0)$ for some closed interval $I$ containing $0$ in its interior.  If for some $s,t >0$, $h_{\gamma}$ is decreasing on $(0,t)$ and increasing on $(-s,0)$, then there is a choice of $\tau'\subset \tau$ determining an attracting neighborhood $ V_{\gamma}(\tau',A)$.  Symmetrically, if $h_{\gamma}$ is increasing on $(0,t)$ and decreasing on $(-s,0)$, then there is a choice of $\tau'\subset \tau$ determining an attracting neighborhood $V_{-\gamma}(\tau',A)$.

\begin{definition}\label{complete attracting} Let $\mathcal F$ be a $C^{\infty,0}$ foliation.  A set of holonomy neighborhoods $V_{\gamma_1}(\tau_1,A_1),...,V_{\gamma_n}(\tau_n,A_n)$ for $\mathcal F$ is {\sl spanning} if each leaf of $\mathcal F$ has nonempty intersection with the interior at least one $V_{\gamma_i}(\tau_i,A_i)$.
\end{definition}

\begin{definition}\label{Vcompatible}
Let $V$ be the union of pairwise disjoint holonomy neighborhoods $V_{\gamma_1}(\tau_1,A_1),...,V_{\gamma_n}(\tau_n,A_n)$ for $\mathcal F$.  A $C^{\infty,0}$ foliation $\mathcal G$ in $M$ is called {\sl $V$-compatible with $\mathcal F$}, (or simply {\sl $V$-compatible} if $\mathcal F$ is clear from context) if each $V_{\gamma_i}(\tau_i,A_i)$ is a holonomy neighborhood for $\mathcal G$, with $V$ spanning for $\mathcal G$ if it is spanning for $\mathcal F$.  
\end{definition}

Fix a set of pairwise disjoint holonomy neighborhoods $V_{\gamma_1}(\tau_1,A_1),...,$ $ V_{\gamma_n}(\tau_n,A_n)$ for $\mathcal F$, and let $V$ denote their union.  Let $R_i=R_{\gamma_i}(\tau_i,A_i), 1\le i\le n$, and let $R$ denote the union of the $R_i$.  For each $i, 1\le i\le n$, fix a smooth open neighborhood $N_{R_i}$ of $R_i$ in $V_i$.  Choose each $N_{R_i}$ small enough so that its closure, $\overline{N_{R_i}}$, is a closed regular neighborhood of $R_i$.  Let $N_R$ denote the union of the $N_{R_i}$.

Now, given $V$, $R$ and $N_R$, we further constrain the set of foliations $\mathcal F$ (that we need to approximate by smooth contact structures) to $C^{\infty,0}$ foliations which are smooth on $N_R$.  The following lemma establishes that we can do this with no loss of generality.

\begin{lemma}\label{smoothaboutR}\cite{KR4} Let $\mathcal F$ be a $C^{\infty,0}$ foliation and let $\Phi$ be a smooth flow transverse to $\mathcal F$.  Let $V$ denote the union of a set of pairwise disjoint holonomy neighborhoods for $\mathcal F$ and fix $N_R$ as above.  There is an isotopy of $M$ taking $\mathcal F$ to a $C^{\infty,0}$ foliation which is both $C^0$ close to $\mathcal F$ and smooth on $N_R$.  This isotopy may be taken to preserve both $V$ and the flow lines of $\Phi$ setwise.
\end{lemma}

Next we describe a preferred product parametrization on a closed set containing $V$.  In this paper, we express $S^1$ as the quotient $S^1=[-1,1]/\sim$, where $\sim$ is the equivalence relation on $[-1,1]$ which identifies $-1$ and $1$.

\begin{lemma}\label{prodnbddefn} Let $\mathcal F$ be a $C^{\infty,0}$ foliation and let $\Phi$ be a smooth flow transverse to $\mathcal F$.  Let $V$ denote the union of pairwise disjoint holonomy neighborhoods $V_i=V_{\gamma_i}(\tau_i,A_i), 1\le i\le n,$ for $\mathcal F$, and fix $N_R$ as above.  Suppose $\mathcal F$ is smooth on $N_R$.  Then for each $i, 1\le i\le n$, there is a pairwise disjoint collection of closed solid tori $P_i$ such that $P_i$ contains $V_i$ and there is a diffeomorphism $P_i\to [-1,1]\times S^1\times [-1,1]$ which satisfies the following: 
\begin{enumerate}
\item the flow segments $\Phi\cap P_i$ are identified with the segments $\{(x,y)\}\times [-1,1]$,
\item $A_i$ is identified with $[-1,1]\times S^1\times \{0\}$,
\item $\gamma_i$ is identified with $\{0\}\times S^1\times \{0\}$, 
\item $R_i$ is identified with $[-1,1]\times \{1\sim -1\}\times [-1,1]$, and
\item the restriction of the diffeomorphism to $N_{R_i}$ maps leaves of $\mathcal F$ to horizontal level surfaces $D_z\times \{z\}$, where $D_z=[-1,1]\times (( [1/2,1] \cup [-1,-1/2])/\sim)$.
\end{enumerate}
\end{lemma}

\begin{proof}
Since $V_i$ is homeomorphic to a solid torus, it is contained in a solid torus which is diffeomorphic to $[-1,1]\times S^1\times [-1,1]$, where the diffeomorphism can be chosen to identify $A$ with $[-1,1]\times S^1\times \{0\}$ and the flow segments $\Phi\cap P$ with the segments $\{(x,y)\}\times [-1,1]$.  Moreover, since the restriction of $\mathcal F$ to $V_i \cap N_R$ is a smooth product foliation transverse to vertical fibers, and there is a unique such up to diffeomorphism, this diffeomorphism $P_i \to [-1,1]\times S^1\times [-1,1]$ can also be chosen so that the restriction of the diffeomorphism to $N_R$ maps leaves of $\mathcal F\cap N_R$ to horizontal level surfaces $D_z\times \{z\}$, where $D_z=[-1,1]\times (([1/2,1] \cup [-1,-1/2])/\sim)$.
\end{proof}

\begin{definition} Fix $V$ and $N_R$ as above.  Let $P_i$ and $P_i\to [-1,1]\times S^1\times [-1,1]$ be as given in Lemma~\ref{prodnbddefn}.  Abuse notation and use the diffeomorphism to identify $P_i$ with $[-1,1]\times S^1\times [-1,1]$.  Let $\mathcal P_i$ be the product foliation of $P_i$ with leaves $([-1,1]\times S^1)\times \{t\}$, and call such a foliated solid torus, $(P_i,\mathcal P_i)$, a {\sl product neighborhood of $(V_i;N_{R_i})$.} Letting $P$ denote the union of the $P_i$ and $\mathcal P$ denote the union of the $\mathcal P_i$, call $(P,\mathcal P)$ a {\sl product neighborhood of $(V;N_R)$.} 
\end{definition}

\begin{definition}\label{stronglycompatible}
Let $\mathcal F$ be a $C^{\infty,0}$ foliation and $V$ the union of pairwise disjoint, holonomy neighborhoods $V_{\gamma_i}(\tau_i,A_i), 1\le i\le k,$ for $\mathcal F$.  Let $R$ denote the union of the $R_{\gamma_i}(\tau_i,A_i), 1\le i\le k$, and let $N_R$ be an open regular neighborhood of $R$ in $V$.  Let $(P,\mathcal P)$ be a product neighborhood of $(V;N_R)$.  The foliation $\mathcal F$ is {\sl strongly $(V,P)$ compatible} if
\begin{enumerate}
\item $\mathcal F\cap N_R=\mathcal P\cap N_R$, and 
\item in the coordinates inherited from $P$, $\mathcal F\cap V$ is a product foliation $[-1,1]\times \mathcal F_0$, where $\mathcal F_0$ is a $C^{\infty,0}$ foliation of $V\cap (\{1\}\times S^1\times [-1,1])$ (i.e., $\mathcal F\cap V$ is $x$-invariant).
\end{enumerate}
\end{definition}

Given $V$, $R$ and $N_R$, we will further constrain the set of foliations $\mathcal F$ (that we need to approximate by smooth contact structures) to $C^{\infty,0}$ foliations which are strongly $(V,P)$ compatible for some choice of product neighborhood $(P,\mathcal P)$.  The following lemma establishes that we can do this with no loss of generality; namely, after a small perturbation of $\mathcal F$, it is possible to rechoose the diffeomorphisms $P_i\to [-1,1]\times S^1\times [-1,1]$ so that $\mathcal F=\mathcal P$ on $\overline{N_R}$ and $\mathcal F$ is invariant under translation in the first coordinate.

\begin{lemma}\label{stronglycompatibleF}\cite{KR4} Let $\mathcal F$ be a $C^{\infty,0}$ foliation and let $\Phi$ be a smooth flow transverse to $\mathcal F$.
Let $V$ denote the union of a set of pairwise disjoint holonomy neighborhoods for $\mathcal F$ and fix $N_R$ as above.  There is an isotopy of $M$ which takes $\mathcal F$ to a $C^{\infty,0}$ foliation which is $C^0$ close to $\mathcal F$ and strongly $(V,P)$ compatible for some choice of product neighborhood $(P,\mathcal P)$ of $(V;N_R)$.  This isotopy may be taken to preserve both $V$ and the flow lines of $\Phi$ setwise.
\end{lemma}

For simplicity of exposition, it is sometimes useful to fix a product metric on $P$.  We do so as follows.

\begin{definition}\label{metricpn} Let $V_{\gamma}(\tau,A)$ be a holonomy neighborhood and let $(P,\mathcal P)$ be a product neighborhood of $V_{\gamma}(\tau,A)$.  

Put the standard metric on intervals and let $S^1=[-1,1]/(-1\sim 1)$ inherit its metric from $[-1,1]$.  Then put the product metric on $P= [-1,1]\times S^1\times [-1,1]$.  When $P$ comes equipped with the identification $P=[-1,1]\times S^1\times [-1,1]$ as given in Lemma~\ref{prodnbddefn} together with this product metric, the product neighborhood $(P,\mathcal P)$ is called a {\sl metric product neighborhood of $V_{\gamma}(\tau,A)$.} 
\end{definition}

A metric product neighborhood induces a metric on both $V_{\gamma}(\tau,A)\subset P$ and, via lengths of paths, on $Q_{\gamma}(\tau,A)$.

\begin{prop}\label{deltaproduct2}
Suppose $V_{\gamma}(\tau,A)$ is a holonomy neighborhood with metric product neighborhood $(P,\mathcal P)$.  Choose $\epsilon>0$.  There is a closed transversal $\tau'\subset\tau$ containing $p$ in its interior such that $T(\mathcal F\cap V_{\gamma}(\tau',A))$ is $\epsilon$ $C^0$ close to $T(\mathcal P \cap V_{\gamma}(\tau',A))$.
\end{prop}

\begin{proof} At each point $x$ of $\gamma$, choose a foliation chart of the form $U_x \times (-\delta_x,\delta_x)$ such that $U_x$ is a neighborhood of $x$ in $A$ and $\{x\} \times (-\delta_x,\delta_x)$ is transverse to both $\mathcal P$ and $\mathcal F$.  Since both foliations are at least $C^{\infty,0}$ and contain $A$ as leaf, $\delta_x$ may be chosen small enough that $T\mathcal P$ and $T\mathcal F$ are within $\epsilon$ at all points of $ U_x \times (-\delta_x,\delta_x)$.  Passing to a finite cover, there exists a short enough transversal $\tau'\subset\tau$ about $p$ such that that all its leaf preserving translates sweep out the desired neighborhood $V_{\gamma}(\tau',A)$.
\end{proof}

\begin{definition} Suppose $V_{\gamma}(\tau,A)$ is a holonomy neighborhood with product neighborhood $(P,\mathcal P)$.  Choose a metric on $P$.  The restriction of $\mathcal F$ to $V_{\gamma}(\tau,A)$ is {\sl $\epsilon$-horizontal} if $T(\mathcal F\cap V_{\gamma}(\tau',A))$ is $\epsilon$ $C^0$ close to $T(\mathcal P \cap V_{\gamma}(\tau',A))$.
\end{definition}

\begin{definition} Suppose $V_{\gamma}(\tau,A)$ is a holonomy neighborhood for $\mathcal F$ with metric product neighborhood $(P,\mathcal P)$.  Fix $\epsilon>0$.  $V_{\gamma}(\tau,A)$ is {\sl $\epsilon$-flat} if the length of the longest segment of the restriction of $\Phi$ to $V_{\gamma}(\tau,A)$ is bounded above by $\epsilon$.
\end{definition}

In other words, $V_{\gamma}(\tau,A)$ is $\epsilon$-flat if its maximum ``height" is small relative to both its minimum ``width'' and its minimum ``length".  The next lemma follows immediately from the continuity of $T\mathcal F$.

\begin{lemma} \label{flatisenough1}
Fix $\epsilon>0$.  There is a $\delta>0$ so that if $V_{\gamma}(\tau,A)$ is $\delta$-flat, then the restriction of $\mathcal F$ to $V_{\gamma}(\tau,A)$ is $\epsilon$-horizontal.  \qed
\end{lemma}

\section{Minimal sets}\label{minimal}

Holonomy neighborhoods which are attracting play a particularly important role in the construction of approximating contact structures.  Results from the last section will eventually be used to make $C^0$ approximations in a single (attracting) holonomy neighborhood.  Results in this section will be used to choose or create enough attracting holonomy neighborhoods so that these approximations can be extended to all of $M$.  This is accomplished through the use of minimal sets.

\begin{definition} Let $\mathcal F$ be a $C^0$ foliation.  A {\sl minimal set} of $\mathcal F$ is a non-empty closed subset of $M$ that is a union of leaves and is minimal with respect to inclusion among such sets.
\end{definition}

Equivalently, a subset $\Lambda$ of $M$ is a minimal set if and only if both $\Lambda =\overline{L}$ for some leaf $L$ and $\Lambda$ properly contains no leaf closure.  In particular, if $\Lambda$ is a minimal set, then $\Lambda=\overline{L}$ for all leaves $L$ contained in $\Lambda$.  More generally, if $L$ is any leaf of $\mathcal F$, then by Zorn's Lemma, $\overline{L}$ contains at least one minimal set.

The next two results follow from elementary point set topology.  See, for example, Theorem~1, Chapter~3, of \cite{goodman}.

\begin{prop}
Let $\Lambda$ be a minimal set of $\mathcal F$.  Then every leaf of $\mathcal F$ contained in $\Lambda$ is dense in $\Lambda$.
\qed 
\end{prop}

\begin{prop}\label{transversals}
Let $\Lambda$ be a minimal set of $\mathcal F$ and let $\tau$ be a transversal to $\mathcal F$ satisfying $\Lambda\cap \text{int} \,\tau\ne\emptyset$.  Then the following properties hold.
\begin{enumerate}
\item The intersection $\Lambda\cap \tau$ is either discrete, all of $\tau$, or a Cantor set.

\item $\Lambda\cap \tau=\tau$ if and only if $\mathcal F=\Lambda$ is minimal.

\item $\Lambda\cap \tau$ is discrete if and only if $\Lambda$ consists of a single compact leaf.

\end{enumerate}
Moreover, if $\sigma$ is any other transversal to $\mathcal F$ satisfying $\Lambda\cap \text{int} \,\sigma\ne\emptyset$, then both $\Lambda\cap \tau$ and $\Lambda\cap\sigma$ are of the same type; namely, discrete, all of the transverse interval, or Cantor.  \qed
\end{prop}

A minimal set $\Lambda$ is called {\sl exceptional} if for some transversal $\tau$, $\Lambda\cap\tau$ is a Cantor set.  A foliation is called {\sl minimal} if it is itself a minimal set.

\begin{definition}
Suppose that $\Lambda$ is a minimal set and that $\tau$ is a transversal to $\mathcal F$.  When $\Lambda\cap \tau$ is a Cantor set, we differentiate between those points of $\Lambda\cap\tau$ which are endpoints of ``removed intervals'' and those which are not.  The leaves of $\Lambda$ which correspond to points which are endpoints of ``removed intervals'' are isolated on one side in $\Lambda$ and so we refer to them as {\sl boundary leaves}.  The leaves $L$ of $\Lambda$ which do not correspond to endpoints of ``removed intervals'' are approached arbitrarily closely on both sides by leaves in $\Lambda$ (and hence by leaves of $L$).  Call these leaves {\sl non-boundary leaves}.  

When $\Lambda\cap \tau=\tau$, Proposition~\ref{transversals} implies $\Lambda = M$ and all leaves of $\Lambda$ are non-boundary leaves.  When $\Lambda\cap \tau$ is discrete, the single compact leaf of $\Lambda$ is a boundary leaf.
\end{definition}

\begin{definition} Let $p\in M$.  A closed segment $\tau$ is a {\sl transversal about $p$} if $\tau$ is a transversal to $\mathcal F$ which contains $p$ in its interior.
\end{definition}

\begin{lemma} \label{fiberthrough}
Let $S\times\{0,1\}$ be leaves of $\mathcal F$ which bound an $I$-bundle $S\times I$, such that each $I$-fiber is transverse to $\mathcal F$.  If $\Lambda\subset S\times I$ is a minimal set of $\mathcal F$, then necessarily $S$ is compact and $\Lambda$ is isotopic to $S\times \{0\}$.  
\end{lemma}

\begin{proof}

Begin by noting that since any $I$-fiber $\{x\}\times (0,1)$, $x\in S$, is disjoint from $S\times\{0,1\}$, 
$\overline{S\times\{0\}}$ has empty intersection with $S\times (0,1)$.  Hence, either $S$ is compact and $\Lambda$ is a component of $S\times\{0,1\}$, or $S\times\{0,1\}$ is disjoint from $\Lambda$.  

Restrict attention therefore to the case that $S\times\{0,1\}$ is disjoint from $\Lambda$.  Since $\Lambda$ is closed, for each $x\in S$ there is a minimum point $m_x$ of $(\{x\} \times I) \cap \Lambda$.  Let $L_\text{min}$ be collection of all such points $m_x$ as $x$ ranges over $S$.  In the analogous way define $L_\text{max}$ to be the union of points $M_x$.

Notice that if $B=D\times I\subset S\times I$ is any flow box, each of $L_\text{min}\cap B$ and $L_\text{max}\cap B$ is a component of $\mathcal F\cap B$.  It follows that $L_\text{min}$ and $L_\text{max}$ are leaves of $\Lambda$.  If for some $x, m_x < M_X$ it follows that $\overline L_\text{min} \cap \overline L_\text{max} = \emptyset$, a contradiction.  Thus $L_\text{min} = L_\text{max} = \Lambda$ is compact, and isotopic to $S\times\{0\}$.  It follows that $S$ is necessarily compact.
\end{proof}

In Section~\ref{noholo}, we will investigate foliations with only trivial holonomy, and the following lemma will prove useful.  It is included here as its proof uses flow box decompositions but is not needed for any of the results of this section.

\begin{lemma} \label{noholIbundle} 
Let $S\times\{0,1\}$ be leaves of $\mathcal F$ which bound an $I$-bundle $S\times I$ in $M$.  Suppose that $\mathcal F$ lies transverse to the $I$-fibers and has only trivial holonomy.  Then the restriction of $\mathcal F$ to $S\times I$ is, up to a $\Phi$-preserving isotopy, the product foliation $S\times I$.
\end{lemma}

\begin{proof}
Let $L$ be any leaf of the restriction of $\mathcal F$ to $S\times (0,1)$.  We claim that $L$ intersects each flow box $D\times I$ in exactly one component, and hence is isotopic, via a $\Phi$ preserving isotopy, to $L\times \{0\}$.  Consider a component $\Delta=D\times \{a\}$ of the intersection of $L$ with flow boxes $B=D\times I$.  Suppose by way of contradiction, that $L\cap B$ contains a second component $\Delta'=D\times \{a'\}$.  

Let $x\in D$.  Since $L$ is connected, there is a path $\rho$ in $L$ from $(x,a)$ to $(x,a')$.  Express $\rho$ as a concatenation of finitely many intervals, each of which lies in a single flow box, and consider the immersed cylinder $\rho\times I$.  A standard cut and paste argument reveals that we may assume that this cylinder is embedded.  But the existence of the path $\rho$ implies the existence of nontrivial holonomy, a contradiction.  Hence, the intersection of $L$ with any flow box is connected.  In addition, since $S$ is path connected, and therefore any two flow boxes are connected by a path in $S\times [0,1]$, the intersection of $L$ with any flow box is also nonempty.  Hence, the restriction of $\mathcal F$ to $S\times [0,1]$ is, up to a $\Phi$-preserving isotopy, the product foliation $S\times I$.
\end{proof}

\begin{cor}\label{finiteness} Suppose $\mathcal F$ is a $C^0$ foliation of a compact 3-manifold $M$.
\begin{enumerate}

\item At most finitely many isotopy classes of compact surface can be realized as leaves of $\mathcal F$.

\item $\mathcal F$ can contain only finitely many exceptional minimal sets.  

\end{enumerate}
\end{cor}

\begin{proof} If $\mathcal F$ is not transversely oriented, work instead with $(\hat{M},\hat{\mathcal F})$, where $\pi:\hat{M}\to M$ is a double cover such that $\hat{\mathcal F}=\pi^{-1}(\mathcal F)$ is transversely oriented.  (See, for example, Proposition~3.5.1 of \cite{CC}.) In this case, $\Lambda$ is a minimal set for $\mathcal F$ only if $\pi^{-1}(\Lambda)$ contains a minimal set for $\hat{\mathcal F}$, and hence if the claimed result holds true for $\hat{\mathcal F}$, it holds true for $\mathcal F$ also.

So restrict attention to the case that $\mathcal F$ is transversely oriented, and let $\Phi$ be an oriented flow transverse to $\mathcal F$.  Choose a flow box decomposition $\mathcal B$ for $(M, \mathcal F, \Phi)$.  Denote a flow box in $\mathcal B$ by $B=D\times I$.  For each leaf $L$ of $\mathcal F$ and flow box $B$ of $\mathcal B$, call any component of $L\cap B$ a {\sl plaque} of $L$.

Let $\mathcal L = \{\Lambda_{\alpha}\}$ denote any finite set of distinct minimal sets of $\mathcal F$, and let $\Lambda=\cup_{\alpha} \Lambda_{\alpha}$.  Let $B_i$ for $i=1, \dots, n$ be those flow boxes such that $B_i \cap \Lambda \ne \emptyset$.  Let $Y$ denote the open manifold $M\setminus \Lambda$.  Let $\mathcal S=\{D_i\times \{0\}| 1\le i\le n\}\cup \{D_i\times \{1\}| 1\le i\le n\}.$ Define a map $f:\mathcal S\to \mathcal L$ by $f(D_i\times \{j\})=\Lambda_{\alpha_0}$ if the plaque of $\Lambda$ in $B_i$ closest to $D_i\times \{j\}$ lies in a leaf of $\Lambda_{\alpha_0}$.  Since $\mathcal S$ is finite, so is the image of $f$.

Consider a minimal set $\Lambda_{\beta}$ in $\mathcal L$ which is not in the image of $f$.  We claim that $\Lambda_{\beta}$ is a compact leaf $L_{\beta}$ and lies as a section of a trivial $\mathbb R$-bundle component of $Y\cup \Lambda_{\beta}$, where the $\mathbb R$-fibers are subsegments of the flow $\Phi$.  

To see that this is true, proceed as follows.  Let $\Sigma$ denote the component of the open manifold $Y\cup \Lambda_{\beta}=M\setminus \cup_{\alpha\ne\beta} \Lambda_{\alpha}$ which contains $\Lambda_{\beta}$.  Choose any $p\in L_{\beta}$, a leaf of $ \Lambda_{\beta}$.  This point $p$ lies in a plaque $D_i\times \{b\}$ of $\Lambda_{\beta}$, for some $0<b<1$.  Consider the two plaques of the closed set $\Lambda\setminus\Lambda_{\beta}$ which lie closest to $D_i\times \{b\}$.  They cobound a flow box $\Sigma_i(a,c) = D_i\times (a,c)\subset\Sigma$, for some $0\le a<b<c\le 1$.  Let $L_a$ be the leaf containing $D_i\times \{a\}$, and let $L_c$ be the leaf containing $D_i\times \{c\}$.

Next consider any $B_j\in \mathcal B$ satisfying $\partial_v B_j\cap\Sigma_i(a,c)\ne\emptyset$.  Since $\Lambda_{\beta}$ is not in the image of $f$, $\Sigma_i(a,c)\cap\partial_h B_j=\emptyset$, and hence $\Sigma_i(a,c)$ naturally extends to an $I$-bundle of the form $(D_i\cup D_j)\times \mathbb R\subset\Sigma$.  Repeat this process of moving through adjacent flow boxes to obtain an exhaustion of $L_a$ by an increasing union of plaques, and an exhaustion of $\Sigma$ by an increasing union of $\mathbb R$-bundles, where the $\mathbb R$-fibers are subsegments of the flow $\Phi$ and these $\mathbb R$-bundles are trivial since the flow is oriented.  To see that all of $\Sigma$ is realized by this exhaustion, note that $\Sigma$ is path connected, and any closed path can be decomposed as a piecewise union of finitely many intervals, each of which lies in a single flow box.  It follows that $\Sigma$ is an a trivial $\mathbb R$-bundle over $L_a$, with all $\mathbb R$-fibers subsegments of the flow $\Phi$.  So the leaves $L_a$, $L_{\beta}$ and $L_c$ are isotopic.  Since $L_a$ and $L_{\beta}$ bound an $I$-bundle and $\overline{L_a}\ne\overline{L_{\beta}}$, it follows from Lemma~\ref{fiberthrough} that $L_a$, $L_{\beta}$ and $L_c$ are compact.

Hence, any exceptional minimal set lies in the image of $f$, and thus, there can be at most $2|\mathcal B|$ exceptional minimal sets.  Moreover, any compact leaf has an isotopy representative that lies in the image of $f$, and thus, at most finitely many isotopy classes of compact surface can be realized as leaves of $\mathcal F$.
\end{proof}

\begin{remark} Recall that any compact 3-manifold admits a triangulation \cite{Mo}.  Moreover, given a foliation $\mathcal F$, there is a triangulation $\mathcal T$ of $M$ which is compatible with $\mathcal F$; namely, $\mathcal F$ is in Haken normal form with respect to $\mathcal T$ (Lemma~1.3, \cite{Ga2}).  Hence there is a similar proof of Corollary~\ref{finiteness} which uses triangulations instead of flow box decompositions.  Our proof doesn't depend on first finding a good triangulation.
\end{remark}

\begin{cor}\label{getaway}
Suppose $\Lambda$ is a minimal set and $p\in\Lambda$.  Let $\tau$ be any transversal about $p$.  If $\Lambda$ is a compact leaf, then there is a transversal $\tau_0\subset\tau$ about $p$ such that if a minimal set intersects $\tau_0$, it is a compact leaf isotopic to $\Lambda$.  If $\Lambda$ is not a compact leaf, there is a transversal $\tau_0\subset\tau$ about $p$ which does not intersect any minimal set other than $\Lambda$.  
\end{cor}

\begin{proof} Since minimal sets are closed and there are only finitely many exceptional minimal sets by Corollary~\ref{finiteness}, it is possible to choose $\tau_0$ about $p$ that misses all exceptional minimal sets that are not equal to $\Lambda$.  Furthermore, Corollary~\ref{finiteness} implies that if there is no transversal $\tau_0$ about $p$ that is disjoint from all minimal sets not equal to $\Lambda$, there must exist a sequence of compact leaves $L_i$ all mutually isotopic that limit on $\Lambda$.  However, referring to the proof of Corollary~\ref{finiteness}, at most finitely many of the compact leaves $L_i$ lie in the image of $f$ and any other, $S=L_j$ say, lies as a surface fiber in one of finitely many $I$-bundles of the form $S\times I$, where the surfaces $S\times \{0,1\}$ are leaves of $\mathcal F$.  Hence, the nonexistence of $\tau_0$ implies that $\Lambda$ is a minimal set embedded in an $I$-bundle $S\times [0,1]$ transverse to the $I$-fibers, and therefore by Lemma~\ref{fiberthrough} is a compact leaf isotopic to $S=L_j$.
\end{proof}

\begin{notation}\label{minsetlisting}
Fix a $C^0$ foliation $\mathcal F$.  If $\mathcal F$ is minimal, set $\Lambda_1=\mathcal F$.  Otherwise, let $\Lambda_1,...,\Lambda_r$ denote the exceptional minimal sets of $\mathcal F$ and let $[L_1],...,[L_s]$ denote the isotopy classes of compact leaves of $\mathcal F$.

For each compact leaf of $\mathcal F$, let $X(L)$ denote the minimal $\mathcal F$-saturated closed subset of $M$ containing all leaves of $\mathcal F$ which are isotopic to $L$.  
\end{notation}

\begin{cor}\label{finitespan}
Let $\mathcal F$ be a $C^0$ foliation that is not a fibering of $M$ over $S^1$.  Then there are finitely many pairwise disjoint holonomy neighborhoods $V_{\gamma_1}(\tau_1,A_1),...,V_{\gamma_{r+s}}(\tau_{r+s},A_{r+s})$ such 
that 
\begin{enumerate}
\item for $1\le i\le r$, $\gamma_i$ is an essential simple closed curve in a leaf of $\Lambda_i$,
\item for $1\le i\le s$, $\gamma_{r+i}$ is an essential simple closed curve in a leaf in the isotopy class $[L_i]$, and each leaf isotopic to $L_i$ lies in the interior of $V_{\gamma_{r+i}}(\tau_{r+i}, A_{r+i})$, and
\item each minimal set of $\mathcal F$ has nonempty intersection with the interior of exactly one $V_{\gamma_i}(\tau_i,A_i)$.
\end{enumerate}
In particular, the set of holonomy neighborhoods $\{V_{\gamma_i}(\tau_i,A_i)\}_i$ is spanning.
\end{cor}

Conditions (1) and (2) guarantee that each minimal set of $\mathcal F$ has nonempty intersection with the interior of at least one $V_{\gamma_i}(\tau_i,A_i)$.  So condition (3) guarantees that each $V_{\gamma_i}(\tau_i,A_i)$, $1\le i\le r$, has nonempty intersection with exactly one minimal set, $\Lambda_i$, and each $V_{\gamma_{r+i}}(\tau_{r+i},A_{r+i})$, $1\le i\le s$, has nonempty intersection with exactly one isotopy class of minimal set, $[L_i]$.

\begin{proof} For each $i$, $1\le i\le r$, let $\gamma_i$ be an essential simple closed curve in a leaf of $\Lambda_i$, and choose a holonomy neighborhood $V_{\gamma_i}(\tau_i,A_i)$.  Choose the $(\tau_i,A_i)$ so that the neighborhoods $V_{\gamma_i}(\tau_i,A_i)$ are pairwise disjoint and disjoint from any compact leaf of $\mathcal F$.  This is possible by Proposition~\ref{getaway}.

Let $L$ be a compact leaf of $\mathcal F$.  Recall that $X(L)$ denotes the minimal $\mathcal F$-saturated closed subset of $M$ containing all leaves of $\mathcal F$ which are isotopic to $L$.  Since $\mathcal F$ is not a fibering of $M$ over $S^1$, either $X(L)= L$ or $X(L)\cong L\times [0,1]$, with the identification given by a diffeomorphism.  Note that if $L $ and $F$ are compact leaves of $\mathcal F$, then either $L$ and $F$ are isotopic and $X(L)=X(F)$, or $L$ and $F$ are not isotopic and $X(L)\cap X(F)=\emptyset$.  Rechoose the isotopy class representatives as necessary so that $X(L_1),...,X(L_s)$ is a listing of the sets $X(L)$, where either $X(L_i)=L_i$ or, under the identification $X(L_i)=L_i\times I$, $L_i$ is identified with $L_i\times \{0\}$.  

Set $n=r+s$.

For each $j$, let $\gamma_{r+j}$ be an essential simple closed curve in $L_j = L_j\times\{0\}$.  Choose a holonomy neighborhood $V_{\gamma_{r+j}}(\tau_{r+j},A_{r+j})$.  If $X(L_j)=L_j\times [0,1]$, choose the transversal $\tau_{r+j}$ long enough so that its interior has nonempty intersection with $L_j\times \{1\}$ (and hence also with all leaves isotopic to $L_j$), and short enough so that the holonomy neighborhoods $V_{\gamma_{r+1}}(\tau_1,A_1),...,V_{\gamma_n}(\tau_n,A_n)$
are pairwise disjoint and disjoint from the minimal sets $\Lambda_1,...,\Lambda_r$ and their fixed holonomy neighborhoods $V_{\gamma_1}(\tau_1,A_1),...,V_{\gamma_r}(\tau_r,A_r)$.

Since the closure of any leaf contains a minimal set, each leaf of $\mathcal F$ has nonempty intersection with the interior of $V_{\gamma_i}(\tau_i,A_i)$ for some $i$; in other words, the collection of holonomy neighborhoods $\{V_{\gamma_i}(\tau_i,A_i)\}_i$ is spanning.
\end{proof}

\section {Generalized Denjoy blow up}\label{Denjoy}
In this section we define the operation of generalized Denjoy blow up.  Informally, this operation consists of thickening a leaf and inserting a new foliation into the thickened region.  This will be used to modify a foliation and create attracting neighborhoods of curves in a leaf of a foliation.

\begin{definition}\label{blowupdef} Let $L$ be a finite (or even countably infinite) union of leaves of a $C^{k,0}$ foliation $\mathcal F$ of $M$ with $k\ge 1$, and let $\Phi$ be a smooth flow transverse to $\mathcal F$.  A $C^{k,0}$ foliation, $\mathcal F'$, is a {\sl generalized Denjoy blow up} of $\mathcal F$ along $L$ if there is an open subset $U\subset M$ and a collapsing map $h:M \to M$ satisfying the following properties:
\begin{enumerate}

\item $\mathcal F'$ is transverse to $\Phi$,

\item there is an injective $C^k$ immersion $j:L \times (0,1) \to M$ such that $j (L\times (0,1)) = U$,

\item for each $x\in L$, $j(\{x\} \times I)$ is contained in a flow line of $\Phi$,

\item $j(L \times \{0\})$ and $j(L \times \{1\})$ are leaves of $\mathcal F'$,

\item $h^{-1}(x)$ is a point for $x \notin L$ and equals $j(\{x\} \times I)$ for $x \in L$,

\item $h$ preserves flow lines of $\Phi$ and maps leaves of $\mathcal F'$ to leaves of $\mathcal F$, and

\item $h$ is $C^0$ on $M$ and $C^k$ when restricted to any leaf of $\mathcal F'$.
\end{enumerate}
When the restriction of $\mathcal F'$ to $j(L\times [0,1])$ is a product foliation, $\mathcal F'$ is also referred to as a {\sl Denjoy blow up} of $\mathcal F$ along $L$.  
\end{definition}

\begin{thm}\cite{Di}\cite{KR4}\label{blowupthm} Let $\mathcal F$ be $C^{k,0}$ foliation with $k\ge 1$ that is transverse to a smooth flow $\Phi$.  Let $L$ be a finite or countable collection of leaves of $\mathcal F$, and let $\mathcal F_1$ be a $C^{k,0}$ foliation of $L \times I$ transverse to the $I$ coordinate that contains $L \times \partial I$ as leaves.  Then there exists $\mathcal F'$ arbitrarily $C^0$ close to $\mathcal F$ that is a generalized Denjoy blow up of $\mathcal F$ along $L$, and such that the pullback of $\mathcal F'$ to $L \times I$ is $C^{k,0}$ equivalent to $\mathcal F_1$.  

Moreover, if $V$ is the union of a set of pairwise disjoint holonomy neighborhoods for $\mathcal F$, $(P,\mathcal P)$ is a product neighborhood of $V$, and $\mathcal F$ is strongly $(V,P)$ compatible, then $\mathcal F'$ can be chosen to be both $V$-compatible with $\mathcal F$ and strongly $(V,P)$ compatible.  
\end{thm}

\begin{remark} The main ideas of Theorem~\ref{blowupthm} are due to Dippolito, \cite{Di}.  We require slightly more than is easily extracted from his work, namely, $C^0$ approximation.  We give a proof of the theorem in \cite{KR4} using flow box decompositions that also allows for $C^1$, rather than $C^{\infty}$, leaves.  We also realize the additional conditions that $\mathcal F'$ be both $V$-compatible with $\mathcal F$ and strongly $(V,P)$ compatible.
\end{remark}

In order to create attracting neighborhoods, we will be interested in inserting foliations into $L\times (0,1)$ of the form described in the following lemma.

\begin{lemma} \label{introspiral} Let $\gamma$ be an oriented essential simple closed curve in $L = L \times \{0\}$.  There is a $C^{\infty,0}$ foliation on $L \times I$, transverse to the $I$-fibers $\{x\}\times I, x\in L$, and such that the holonomy $h_\gamma$ along $\gamma$ is monotone decreasing on the interior of $I$.  Moreover, this foliation of $L \times I$ can be chosen so that $L\times \{0\}$ and $L\times \{1\}$ are its only minimal sets.
\end{lemma}

\begin{proof} If $\gamma$ does not separate, let $\alpha$ be a properly embedded curve in $L$ that intersects $\gamma$ in a point.  Let $\mathcal P$ be the product foliation on $(L \backslash \alpha) \times I$, and let $h_\gamma:I \to I$ be the desired holonomy about $\gamma$.  Then glueing the leaves of $\mathcal P$ at height $x$ to those of height $h_\gamma(x)$ gives a foliation on $L\times I$ with holonomy $h_\gamma$ around $\gamma$.

If $\gamma$ separates $L$ into components $A$ and $B$, first consider the case that $A$ is compact with genus $g\ge 1$.  In the usual way, $A$ may be thought of a disk $D$ with $4g$ disjoint subarcs glued in pairs.  Let $\mathcal P$ be the product foliation on $D \times I$.  Using $2g$ homeomorphisms $h_i$ of $I$ to glue up the leaves of $\mathcal P$ produces a foliation on $A \times I$.  With the right choice of pairings of glued subarcs, the holonomy along $\gamma$ will be the product of $g$ commutators of the $h_i$.  Since any orientation preserving homeomorphism of $I$ can be written as a single commutator, see for instance, Lemma~3.1 of \cite{Li}, this construction can be carried out for any genus and any choice of $h_\gamma$.

Next, consider the case that $A$ is not compact.  Let $\alpha$ be a properly embedded half-infinite line contained in $A$ and starting on $\gamma$.  Splitting $A$ along $\alpha$ and glueing leaves with $h_\gamma$ as in the non-separating case gives the desired foliation around $\gamma$.

The same constructions are used to extend a given choice of holonomy across $B$.
\end{proof}

\begin{notation} Any foliation given by a generalized Denjoy blow up of $\mathcal F$ along $L$, with the foliation inserted into $L\times (0,1)$ of the type generated by Lemma~\ref{introspiral}, is denoted $\mathcal F' = \mathcal F(L, \gamma)$.
\end{notation}

The following two lemmas use the notation of Definition~\ref{blowupdef} to describe the effect of generalized Denjoy blow up on the set of minimal sets of F.

\begin{lemma}\label{nonmimimalblowup} Let $L$ be a leaf of $\mathcal F$ that is not contained in any minimal set of $\mathcal F$.  Let $\mathcal F'$ be a generalized Denjoy blow up of $\mathcal F$ along $L$.  There is a bijective correspondence between the minimal sets of $\mathcal F$ and those of $\mathcal F'$ given by $\Lambda_\beta \longleftrightarrow h^{-1}(\Lambda_{\beta})$.  In particular, neither $L_0,L_1$, nor any leaf in $j(L \times (0,1))$, is contained in a minimal set of $\mathcal F'$.  The restriction of the collapsing function, $h$, gives a homeomorphism from $h^{-1}(\Lambda_\beta)$ to $\Lambda_\beta$ for all $\beta$.
\end{lemma}

\begin{proof} 

Since $\overline{L}$ is not minimal, it properly contains $\overline{S}$ for some leaf $S$ of $\mathcal F$.  It follows that each of $\overline{L}_0$ and $\overline{L}_1$ properly contains $h^{-1}(\overline{S}) = \overline{h^{-1}(S)}$.  So neither $\overline{L}_0$ nor $\overline{L}_1$ is minimal.  For each leaf $F\subset j(L \times (0,1))$ the set $\overline{F}$ properly contains each of $L_0$ and $L_1$ and hence each of $\overline{L}_0$ and $\overline{L}_1$.  Therefore $\overline F$ is not minimal.

Finally, notice that $S'$ is a leaf of $\mathcal F'$ not equal to $L_0, L_1$ or a leaf of $j(L \times (0,1))$ if and only if $S=h(S')$ is a leaf of $\mathcal F$ not equal to $L$.  Since $h^{-1}(\overline{S}) = \overline{h^{-1}(S)}$, the claimed bijective correspondence of minimal sets follows immediately.
\end{proof}

\begin{lemma} \label{minlcreate} Suppose $\overline{L}$ and $\Lambda$ are minimal sets of $\mathcal F$ with $\overline L \neq \Lambda$.  Let $\mathcal F'$ be a generalized Denjoy blow up of $\mathcal F$ along $L$.  Then $h^{-1}(\Lambda)$ is a minimal set of $\mathcal F'$.  Any other minimal set of $\mathcal F'$ arises in one of the following ways.
\begin{enumerate}
\item If $L$ is compact, let $\mathcal F_1$ be the foliation of $L \times I$ that is contained in $\mathcal F'$.  Minimal sets of $\mathcal F_1$ are mapped to minimal sets of $\mathcal F'$ by inclusion; in particular $L_0$ and $L_1$ are minimal sets.

\item If $L$ is noncompact, with $L$ a nonboundary leaf of $\overline{L}$, then $\overline{L}_0=\overline{L}_1$ is a minimal set of $\mathcal F'$.

\item If $L$ is noncompact, with $L$ a boundary leaf of $\overline{L}$, then there are two possibilities, depending on whether $L$ is isolated in $\overline{L}$ from above or from below.  If $L$ is isolated in $\overline{L}$ from below (above), then $\overline{L}_1$ (respectively, $\overline{L}_0)$ is a minimal set, and $\overline{L}_0$ (respectively, $\overline{L}_1)$ properly contains this minimal set.  In particular, the leaf $L_0$ (respectively, $L_1)$ is not contained in a minimal set of $\mathcal F'$.
\end{enumerate}
\end{lemma}

\begin{proof} Since $\mathcal F'$ is not minimal, any minimal set of $\mathcal F'$ is either a compact leaf or else an exceptional minimal set.  Moreover, the only minimal set impacted by the blow up of $L$ is $\overline{L}$.  So we restrict attention to $\overline{L}$ and $h^{-1}(\overline{L})$.

If $L$ is compact, then so is $h^{-1}(\overline{L}) = j(L \times [0,1])$ and (1) follows immediately.  

Suppose instead that $L$ is not compact and $\Lambda'$ is a minimal set of $\mathcal F'$ that intersects $j(L \times (0,1))$.  It follows from Lemma~\ref{fiberthrough} that $\Lambda' \cap (\overline{L}_0\cup\overline{L}_1) \ne \emptyset$, and hence that $\Lambda' \subset (\overline{L}_0\cup\overline{L}_1)$.

Since $h(\overline{L}_0\cup\overline{L}_1) = \overline{L}$, the union $\overline{L}_0\cup\overline{L}_1$ cannot properly contain $\overline{S}$, for some leaf $S$ of $\mathcal F'$, unless $\overline{S}=\overline{L}_i$ for some $ i\in \{0,1\}$.  It follows that $\Lambda'$ is equal to one or both of $\overline{L}_0$ and $\overline{L}_1$.

The cases are distinguished as follows.  If $L$ is a nonboundary leaf, $\overline{L}_0=\overline{L}_1$ is minimal, and if $L$ is a boundary leaf, exactly one of $\overline{L}_0$ and $\overline{L}_1$ is minimal.  To see this, let $\tau$ be a transversal to $\mathcal F$ containing a point $p\in \tau\cap L$.  There is a sequence of points $p_n\in\tau\cap L$ which limit on $p$.  If there exists a sequence of such points limiting on $p$ from above, then the points $h^{-1}(p_n)\cap L_1$ limit on $L_1$ from above.  Similarly, if there exists a sequence of such points limiting on $p$ from below, then the points $h^{-1}(p_n)\cap L_0$ limit on $L_0$ from below.  It follows that $L$ is isolated in $\overline{L}$ from below if and only if $L_0$ is isolated, and that $L$ is isolated in $\overline{L}$ from above if and only if $L_1$ is isolated.  So either $L$ is a nonboundary leaf and $\overline{L}_0=\overline{L}_1$ is minimal, or $L$ is a boundary leaf, isolated from either above or below, and exactly one of $\overline{L}_0$ or $\overline{L}_1$ is minimal.
\end{proof}

\section{Creating attracting holonomy}

In this section, we restrict attention to the case that every minimal set of $\mathcal F$ contains a leaf which is not homeomorphic to $\mathbb R^2$.  The remaining case, in which $M=T^3$ \cite{GICM}, is considered in Section~\ref{noholo}.

Consider a minimal set $\Lambda$ of a foliation $\mathcal F$.  Restrict attention to the case that $\Lambda$ is not a compact leaf.  So either $\Lambda=\mathcal F$ or $\Lambda$ is exceptional.  When the foliation $\mathcal F$ under consideration is $C^2$, a result of Sacksteder~\cite{S} guarantees the existence of a leaf $L$ in $\Lambda$ and simple closed curve $\gamma$ in $L$ such that the holonomy $h_{\gamma}$ along $\gamma$ is linear attracting, that is $h'(0)<1$.  As shown by Eliashberg and Thurston, this combination of smoothness and linear attracting holonomy can be used to introduce a contact region in a neighborhood of $\gamma$.

When the foliation $\mathcal F$ is only $C^{\infty,0}$, it is shown in \cite{KR2} that it is possible to introduce contact regions about a simple closed loop $\gamma$ in a leaf $L$ of $\Lambda$ for which the foliation has a (topologically) attracting neighborhood.  In general, however, there might be no curve with such an attracting region.  

In this section, we show that by taking advantage of generalized Denjoy blow up, it is possible to $C^0$ approximate $\mathcal F$ by a foliation $\mathcal F'$, where each minimal set of $\mathcal F'$ has nonempty intersection with one of finitely many attracting neighborhoods, $V_{\gamma_1}(\tau_1,A_1),...,V_{\gamma_n}(\tau_n,A_n)$.  Moreover, $\mathcal F'$ and the attracting neighborhoods $V_{\gamma_1}(\tau_1,A_1),...,V_{\gamma_n}(\tau_n,A_n)$ can be chosen so that, for each $i$, $\tau_i$ is small and the restriction of $\mathcal F'$ to $V_{\gamma_i}$ is $C^0$ close to a product foliation $A_{\gamma_i}\times\tau_i$.

In order to make sense of ``small'' and $C^0$ close, it is useful to fix a Riemannian metric $g$ on $M$.  We choose a particularly convenient $g$ as follows.  Recall that if $\mathcal F$ is a $C^{\infty,0}$ foliation of $M$ which is not a fibering of $M$ over $S^1$, then Corollary~\ref{finitespan} guarantees the existence of a finite spanning set of pairwise disjoint holonomy neighborhoods for $\mathcal F$.  

\begin{notation} \label{metricmetric} Let $\mathcal F_0$ be a $C^{\infty,0}$ foliation of $M$.  If $\mathcal F_0$ is a fibering of $M$ over $S^1$, perform a $C^0$ small Denjoy splitting of $\mathcal F_0$ along a fiber and let $\mathcal F_1$ denote this new $C^{\infty,0}$ foliation.  If $\mathcal F_0$ is not a fibering of $M$ over $S^1$, let $\mathcal F_1 = \mathcal F_0$.

Let $\{V'_1,...,V'_n\}$, $n=r+s$, denote a set of pairwise disjoint holonomy neighborhoods $V'_i=V_{\gamma_i}(\sigma_i,A_i)$ for $\mathcal F_1$ satisfying the conditions of Corollary~\ref{finitespan}.  Let $V'$ denote the union $V'=\cup_i V'_i$.  For each $i, 1\le i\le n$, let $R'_i=R_{\gamma_i}(\sigma_i,A_i)$, and set $R'=\cup_i R_i$.  For each $i, 1\le i\le n$, fix a smooth open neighborhood $N_{R_i'}$ of $R_i'$ in $V_i'$.  Choose each $N_{R_i'}$ small enough so that its closure, $\overline{N_{R_i'}}$, is a closed regular neighborhood of $R_i'$.  Let $N_R'$ denote the union of the $N_{R_i'}$.

Apply Lemma~\ref{stronglycompatibleF} to isotope $\mathcal F_1$ to a $C^{\infty,0}$ foliation $\mathcal F_2$ which is $C^0$ close to $\mathcal F_1$ and strongly $(V',P)$ compatible for some choice of product neighborhood $(P,\mathcal P)$ of $(V';N_{R'})$.

Put the product metric on each component $P_i=[-1,1]\times S^1\times [-1,1]$ of $P$, as described in Definition~\ref{metricpn}, and let $g_0$ denote the resulting metric on $P$.  Let $g=g(P)$ be any fixed Riemannian metric on $M$ which restricts to $g_0$ on $P$.  Since the metric product neighborhoods $P_i$ have pairwise disjoint closures, a partition of unity argument can be used to construct such a metric $g(P)$.
\end{notation}

Beginning with a minimal set $\Lambda$, a leaf $L$ of $\Lambda$, and a simple closed curve $\gamma$ in $L$ which is not homotopically trivial, we show how to introduce an attracting neighborhood, or, sometimes a pair of attracting neighborhoods, about $\gamma$ via generalized Denjoy blow up.  These operations are performed without increasing the number of minimal sets.  Since the goal is to produce $\epsilon$-flat holonomy neighborhoods, it may be necessary, as in (2) of Theorem~\ref{metriccreateholonbhds} below, to introduce new holonomy neighborhoods to take care of thick collections of parallel compact leaves.  

\begin{thm}\label{metriccreateholonbhds} Let $\mathcal F_0$ be a $C^{\infty,0}$ foliation of $M$.  Let $\gamma_1,...,\gamma_r,...,\gamma_{n=r+s}$, $V'$, $\mathcal F_2$, $(P,\mathcal P)$, and $g=g(P)$ be as given in Notation~\ref{metricmetric}.  So, in particular, $\mathcal F_2$ is not a fibering, $V'$ is spanning for $\mathcal F_2$, and $\mathcal F_2$ is strongly $(V',P)$ compatible.

Fix $\epsilon>0$.  There is a $C^{\infty,0}$ foliation $\mathcal F$ that is $\epsilon$ $C^0$ close to $\mathcal F_2$, $V'$-compatible with $\mathcal F_2$, and strongly $(V',P)$ compatible, and a finite set of pairwise disjoint attracting neighborhoods $$V_{\gamma_1}(\tau_1,A_1),...,V_{\gamma_m}(\tau_m,A_m), $$ $m\ge n$, for $\mathcal F $ such that 
\begin{enumerate}

\item for each $i, 1\le i\le n$, $V_{\gamma_i}(\tau_i,A_i)\subset V'_{\gamma_i}(\sigma_i,A_i)$ 

\item for each $i, n<i\le m,$ $\gamma_i$ lies in a compact leaf of $\mathcal F$ and is isotopic to $\gamma_{j_i}$ for some $r<j_i\le n$, and $V_{\gamma_i}(\tau_i,A_i)\subset V'_{\gamma_{j_i}}(\sigma_{j_i},A_{j_i})$, 

\item each $V_{\gamma_i}(\tau_i,A_i)$ is $\epsilon$-flat with respect to $\mathcal F$, 

\item the restriction of $\mathcal F $ to any $V_{\gamma_i}(\tau_i,A_i)$ is $\epsilon$-horizontal, 

\item there is a regular neighborhood $N_h$ of $\partial_h V$ such that $\overline{N_h}\cap A=\emptyset$ and the restriction of $\mathcal F$ to $N_h$ is a smooth product foliation, and 

\item each minimal set of $\mathcal F $ has nonempty intersection with the interior of exactly one $V_{\gamma_i}(\tau_i,A_i)$.  
\end{enumerate}
\end{thm}

Notice that since $\mathcal F$ is strongly $(V',P)$ compatible, it is strongly $(V,P)$ compatible.
Notice also that condition (6) implies that the collection of attracting neighborhoods $V_{\gamma_1}(\tau_1,A_1),...,V_{\gamma_m}(\tau_m,A_m)$ is spanning.  Referring back to Notation~\ref{metricmetric}, each $\gamma_i$ is an essential loop in a leaf of a minimal set.  Therefore each attracting neighborhood $V_{\gamma_i}(\tau_i,A_i)$ necessarily has nonempty intersection with at least one minimal set.  In addition, by the choice of $V'$ and Corollary~\ref{finitespan}, (1) and (2) guarantee that if an attracting set $V_{\gamma_i}(\tau_i,A_i)$, has nonempty intersection with distinct minimal sets, then necessarily, $i>r$ and the minimal sets are isotopic compact leaves.

\begin{proof} By Corollary~\ref{flatisenough1}, it is sufficient to to prove that there is a $C^{\infty,0}$ foliation $\mathcal F$ arbitrarily $C^0$ close to $\mathcal F_2$ such that $\mathcal F$ admits a finite spanning set of pairwise disjoint, $\epsilon$-flat attracting neighborhoods $V_{\gamma_1}(\tau_1,A_1),...,V_{\gamma_n}(\tau_n,A_n)$ for $\mathcal F$ satisfying conditions (1)-(6).  Since the holonomy neighborhoods $V'_{\gamma_1}(\sigma_1,A_1),...,V'_{\gamma_n}(\sigma_n,A_n)$ are pairwise disjoint, condition (1) will guarantee that the neighborhoods $V_{\gamma_1}(\tau_1,A_1),...,$ $V_{\gamma_n}(\tau_n,A_n)$ are pairwise disjoint.

Using Notation~\ref{minsetlisting}, $\mathcal F_2$ has finitely many minimal sets $\Lambda_1,...,\Lambda_r$ that have no compact leaves, and at most finitely many isotopy classes $[L_1],...,[L_s]$, of compact leaves.

Let $L$ be a compact leaf of $\mathcal F_2$.  Since $\mathcal F_2$ is not a fibering over $S^1$, either $X(L)=L$ or $X(L)=L\times [0,1]$, with the identification given by a diffeomorphism.  Abuse notation and set $\Lambda_{r+j}=X(L_j), 1\le j\le s$.

Inductively create a $V'$-compatible, $C^{\infty,0}$ foliation $\mathcal F^k$ arbitrarily $C^0$ close to $\mathcal F_2$ and a pairwise disjoint collection of attracting neighborhoods $$\mathcal V_k=\{V_{\gamma_1}(\tau_1,A_1),...,V_{\gamma_k}(\tau_{n_k},A_{n_k})\}$$ for $\mathcal F^k$ satisfying conditions (1)--(4).  Let $N_0(\mathcal F^k)$ be the number of minimal sets $\Lambda_j, 1\le j\le r,$ that do not intersect the interior of some $V_{\gamma_i}(\tau_i,A_i)$.  Let $N_1(\mathcal F^k)$ be the number of isotopy classes of compact leaves of $\mathcal F^k$ containing leaves that do not intersect the interior of some $V_{\gamma_i}(\tau_i,A_i)$.  Set $N(\mathcal F^k)=N_0(\mathcal F^k)+N_1(\mathcal F^k)$.

The construction of the desired set of attracting neighborhoods consists of finding, or creating, attracting $\epsilon$-flat neighborhoods that decrease $N(\mathcal F^k)$.  These neighborhoods are created by performing, as necessary, generalized Denjoy blow ups very close to leaves $L$ in the sets $\Lambda_i, 1\le i\le r+s$, and by Theorem~\ref{blowupthm}, these blow ups can be chosen arbitrarily $C^0$ close to $\mathcal F$.  The construction takes different forms depending on properties of $L$ and is carried out in Propositions~\ref{compactholonbhd}, \ref{noncompactholonbhd}, and \ref{smallattracting22}.

At the $k^{\text{th}}$ stage in the induction, one or more $\epsilon$-flat attracting neighborhoods are added to $\mathcal V_k$ to yield $\mathcal V_{k+1}$.  
\end{proof}

To complete the proof of Theorem~\ref{metriccreateholonbhds}, it suffices to establish Propositions~\ref{compactholonbhd}, \ref{noncompactholonbhd}, and \ref{smallattracting22}.  We now do so.

Let $\gamma$ be an oriented essential simple closed curve in a leaf $L$ contained in a minimal set $\Lambda$ of $\mathcal F^k$ that contributes to $N(\mathcal F^k)$.  To simplify notation, let $\mathcal F=\mathcal F^k$ at the kth step.  We begin by considering a holonomy neighborhood $V_{\gamma}(\tau,A)\subset P$.  Using the notation of Section~\ref{background}, $p\in \gamma$ and $\tau$ and $\sigma$ are transversals through $p$ such that $h_{\gamma}:\tau\to \sigma$ is a holonomy map for $\mathcal F$ along $\gamma$.  Notice that if $\tau$ is chosen to be sufficiently small, then $V_{\gamma}(\tau,A)$ is $\epsilon$-flat.

\begin{lemma}\label{atnbd} Let $V_{\gamma}(\tau,A)$ be a holonomy neighborhood of $\gamma$.  One of the following is true:
\begin{enumerate}
\item There is a choice of $\tau'\subset \tau$ so that one of $V_{\gamma}(\tau', A)$ and $V_{-\gamma}(\tau', A)$ is an attracting neighborhood, and including the attracting neighborhood decreases $N(\mathcal F)$.

\item There is a choice of $\tau'\subset \tau$ so that, after performing a generalized Denjoy blow up along $L$ arbitrarily $C^0$ close to the identity and compatible with $V_{\gamma}(\tau, A)$, $V_{\gamma}(\tau, A)$ is the union of two attracting neighborhoods $V_{\gamma_0}$ and $V_{\gamma_1}$, where $\gamma_0$ and $\gamma_1$ are the copies of $\gamma$ obtained by the splitting, chosen with opposite orientations as determined by the form of $V_{\gamma}$.  Including these attracting neighborhoods decreases $N(\mathcal F)$.

\item The holonomy map $h_{\gamma}$ is the identity when restricted to at least one of the components of $\tau\backslash {\{p\}}$.
\end{enumerate}
\end{lemma}

\begin{proof}
Consider the holonomy map $h_{\gamma}:\tau\to \sigma$.  If there are intervals in each component of $\tau\backslash {\{p\}}$ on which $h_{\gamma}$ is strictly monotonic, then (1) or (2) must hold.  Otherwise, (3) holds.
\end{proof}

By the lemma, it is enough to consider the case that $\gamma$ is essential in $L$, and $h_{\gamma}$ is the identity when restricted to at least one of the components of $\tau\backslash {\{p\}}$.  Notice that since $\mathcal F$ is taut, $L$ is $\pi_1$-injective, and hence $\gamma$ is homotopically nontrivial in $M$.  

Identify $(\tau,p)$ with $([-u,v],0)$.  For $t\in [-u,v]$, let $L_t$ denote the leaf of $\mathcal F$ intersecting $\tau$ at height $t$.  Call a compact leaf of $\mathcal F$ {\sl isolated} if there is an $\mathcal F$ saturated open neighborhood of $L$ containing no other compact leaf.  For simplicity of exposition, we will consider first the case that $L$ is either noncompact or else compact but isolated.  

As a first step towards building an attracting neighborhood, we show that without increasing $N(\mathcal F)$, generalized Denjoy blow ups can be used to replace holonomy that is the identity on one side of $\gamma$ with monotone holonomy.

\begin{lemma}\label{compactisolated} Let $\gamma$ be an oriented essential simple closed curve in an isolated compact leaf $L$ such that $h_\gamma$ is the identity on $[0,v]$.  Then there exists a generalized Denjoy blow up of $\mathcal F$ to $\mathcal F'$ such that in $\mathcal F'$, $h_\gamma$ is strictly monotone on a nondegenerate subinterval of $[0,v]$.  Moreover, $N(\mathcal F') \leq N(\mathcal F)$, and $h_\gamma$ may be created to be either monotone increasing or decreasing on this subinterval.
\end{lemma}

\begin{proof}
By Corollary~\ref{getaway}, there exists $w\in (0,v]$ small enough to guarantee that the transversal $(0,w)$ is disjoint from the minimal sets of $\mathcal F$.  Let $t\in (0,w)$, and consider the leaf $L_t$.

Let $\gamma_t$ be a parallel copy of $\gamma$ lying in $L_t$ and passing through $t$.  Perform generalized Denjoy blow up $\mathcal F(L_t,\gamma_t)$, arbitrarily $C^0$ close to the identity and compatible with $V_{\gamma}$, to introduce spiraling about $\gamma_t^{\pm 1}$, the two copies of $\gamma_t$ introduced.  
This introduces a nondegenerate interval in $(0,v)$ on which $h_{\gamma}$ is strictly monotonic.  Moreover, the type of monotonicity, increasing or decreasing, can be chosen.  Since $\overline L_t$ is not minimal, Lemma~\ref{nonmimimalblowup} guarantees that no new minimal sets are introduced.
\end{proof}

\begin{prop}\label{compactholonbhd} Let $\gamma$ be an oriented essential simple closed curve in an isolated compact leaf $L$.  If $h_\gamma$ is the identity on at least one of $[-u,0]$ and $[0,v]$, then there exists a generalized Denjoy blow up $\mathcal F'$ of $\mathcal F$ that creates no new compact leaves, has $N(\mathcal F') \leq N(\mathcal F)$ and for which there exists an attracting holonomy neighborhood containing $L$.
\end{prop}

\begin{proof} Apply Lemma~\ref{compactisolated} once or twice, as necessary, and let $\mathcal F'$ be the result of doing the generalized Denjoy blow up or blow ups required to make $h_{\gamma}$ strictly monotone on nondegenerate subintervals of each of $[-u,0]$ and $[0,v]$.  The monotonicity can be chosen so that $h_{\gamma}$ is either attracting on both subintervals or repelling on both subintervals.  Choose $\tau'\subset\tau$ to be the smallest closed interval containing both subintervals.  Including $V_\gamma(\tau')$ with the collection of holonomy neighborhoods may not create intersections with every leaf istotopic to $L$, that is it may not decrease $N_1(\mathcal F)$, but it keeps $N(\mathcal F') \leq N(\mathcal F)$.
\end{proof}

\begin{lemma}\label{noncompact} 
Let $\gamma$ be an oriented essential simple closed curve in a noncompact leaf $L$ of a minimal set such that $h_\gamma$ is the identity on $[0,v]$.  Then there exists a generalized Denjoy blow up of $\mathcal F$ to $\mathcal F'$ such that in $\mathcal F'$, $h_\gamma$ is strictly monotone on a nondegenerate subinterval of $[0,v]$.  Moreover, $N(\mathcal F') \leq N(\mathcal F)$, and $h_\gamma$ may be created to be either monotone increasing or decreasing on this subinterval.
\end{lemma}

\begin{proof} 
By Corollary~\ref{getaway}, there exists $w\in (0,v]$ small enough to guarantee that the transversal $(0,w)$ is disjoint from the minimal sets of $\mathcal F$.  Hence the proof of Lemma~\ref{compactisolated}, the proof in that case works in this case as well.
\end{proof}

\begin{prop}\label{noncompactholonbhd} Let $\gamma$ be an oriented essential simple closed curve in a noncompact leaf $L$ of a minimal set that contributes to $N(\mathcal F)$.  If $h_\gamma$ is the identity on at least one of $[-u,0]$ and $[0,v]$, then there exists a generalized Denjoy blow up $\mathcal F'$ of $\mathcal F$ that creates no new compact leaves, satisfies $N(\mathcal F') < N(\mathcal F)$, and for which there exists an attracting holonomy neighborhood containing $L$.
\end{prop}

\begin{proof}The proof is similar to the proof of Proposition~\ref{compactholonbhd}, instead requiring one or two applications of Lemma~\ref{noncompact}.  The result is a strict decrease in $N_0(\mathcal F)$, thereby forcing $N(\mathcal F') < N(\mathcal F)$.  
\end{proof}

At this point, attracting neighborhoods have been constructed which intersect every minimal set consisting of noncompact leaves (with at least one non-$\mathbb R^2$ leaf).  There remain minimal sets that consist of a single compact leaf.  Lemma~\ref{compactisolated} shows how to construct a holonomy neighborhood that will contain such a leaf.  The next proposition shows how to deal with possibly infinite families of isotopic compact leaves.

Let $L$ be a compact leaf of $\mathcal F$.  Recall that $X(L)$ denotes the minimal $\mathcal F$ saturated closed submanifold of $M$ containing all leaves isotopic to $L$, and that $\mathcal F$ is not a fibering over $S^1$.  So either $X(L)=L$ or $X(L)$ is diffeomorphic to $L\times [0,1]$.  We now restrict attention to the remaining case, that $X(L)\cong L\times [0,1$.  Notice that leaves of $\mathcal F$ that are contained in $X(L)$ are not required to be homeomorphic to $L$.

\begin{prop}\label{smallattracting22} Let $L$ be a compact leaf such that $X(L) \cong L\times [0,1]$, and let $\gamma$ be an oriented essential simple closed curve in $L$.  There is a $C^{\infty,0}$ foliation $\mathcal F'$ $C^0$ close to $\mathcal F$, such that all surfaces of $\mathcal F'$ isotopic to $L$ are covered by finitely many pairwise disjoint, $\epsilon$-flat, attracting neighborhoods that are disjoint from all minimal sets not isotopic to $L$.  

Moreover, such a foliation $\mathcal F'$ can be obtained from $\mathcal F$ by performing a finite number of generalized Denjoy blow ups along leaves of $\mathcal F$ isotopic to $L$, and the attracting neighborhoods can be chosen of the form $V_{\gamma_t}(\tau_t,A_t)$, where $L_t$ is a leaf of $\mathcal F$ isotopic to $L$ and $(L,A,\gamma)$ is isotopic to $(L_t,A_t,\gamma_t)$.  
\end{prop}

\begin{proof}
Let $A$ be a smooth regular neighborhood of $\gamma$ in $L$ and choose $p\in\gamma$.
Use a diffeomorphism to identify $X(L)$ with $L\times [0,1]$.  Choose the diffeomorphism so that it agrees with the product structure of $P$ restricted to $X(L)$.  In particular, flow lines of $\Phi$ map to the vertical fibers $\{x\}\times [0,1]$.  Let $\tau_0$ denote the vertical fiber $\{p\}\times [0,1]$.  Let $\tau$ denote a closed flow segment of $\Phi$ containing $\tau_0$ in its interior, and such that $\tau$ is disjoint from all minimal sets of $\mathcal F$ not isotopic to $L$.

Let $L_t, t\in \mathcal T,$ be a listing of the leaves of $\mathcal F$ isotopic to $L$, where $L_t\cap\tau = \{(p,t)\}$.  Let $\gamma_t=L_t\cap (\gamma\times [0,1])$ and $A_t=L_t\cap (A\times [0,1])$.  Notice that since lengths can only shorten under projection, the minimum length and width of $A_t$ is bounded below by the minimum length and width of $A$.  For all $t\in\mathcal T$, choose an $\epsilon/2$-flat holonomy neighborhood $V_{\gamma_t}(\tau_t,A_t)$, where $\tau_t\subset \tau$.

Since the union $\bigcup_{t\in\mathcal T} A_t$ is compact (since closed), the interiors of finitely many of the $V_{\gamma_t}(\tau_t,A_t)$ cover $\bigcup_{t\in\mathcal T} A_t$.  Choose $t_1,...,t_k\in\mathcal T$ so that the interiors of the $V_{\gamma_{t_i}}(\tau_{t_i},A_{t_i}), 1\le i\le k,$ cover $\bigcup_{t\in\mathcal T} A_t$.  Set $\tau'=\tau_{t_1}\cup...\cup\tau_{t_k}$.  Choose closed subintervals of the $\tau_{t_i}$ and relabel as necessary so that $\tau'=\tau_{t_1}'\cup...\cup \tau_{t_k}'$, where the interiors of the $\tau_{t_i}'$ are pairwise disjoint and only successive $\tau_{t_i}$'s can have nonempty intersection.  Then the $V_{\gamma_{t_i}}(\tau_{t_i}',A_{t_i})$ cover, and after performing finitely many $C^0$ small blow ups along compact leaves of the form $L_a, a\in \partial \tau_{t_i}',$ and then extending the $\tau_{t_i}'$ slightly as necessary, we have the claimed finite set of pairwise disjoint, $\epsilon$-flat, attracting neighborhoods, disjoint from all minimal sets not isotopic to $L$.  
\end{proof}

By Proposition~\ref{finiteness}, only finitely many applications of Propositions~\ref{compactholonbhd} and \ref{smallattracting22} generate attracting holonomy neighborhoods that intersect all compact leaves.  Proposition~\ref{noncompactholonbhd} takes care of the rest of the cases, and therefore completes the proof of Theorem~\ref{metriccreateholonbhds}.  
{\qed}

\section{Approximation in an attracting neighborhood}\label{approximation}

In this section we show (Theorem~\ref{C0V}) how the tangent plane field of the restriction of a $C^{\infty,0}$ foliation $\mathcal F_0$ to a sufficiently small attracting neighborhood $V_{\gamma}(\tau,A)$ can be $C^0$ approximated by a contact structure.  The idea is that restriction of the foliation $\mathcal F_0\cap V_{\gamma}(\tau,A)$ can be analyzed by cutting $V_{\gamma}(\tau,A)$ open along $R_{\gamma}(\tau,A)$ and considering the resulting foliation on $Q_{\gamma}(\tau,A)$.  Much of the foliation data is then encoded in monodromy maps along the four vertical sides of $Q_{\gamma}(\tau,A)$.

To simplify the analysis, we replace $\mathcal F_0$, as necessary, by a $C^0$ close foliation $\mathcal F$ satisfying the conclusions of Theorem~\ref{metriccreateholonbhds}.  

 A key result is Proposition~\ref{decrexists}, in which we show that the restriction of $\mathcal F$ to $\partial_v Q$ can be approximated by a smooth foliation on $\partial_v Q_{\gamma}(\tau,A)$ that has decreasing monodromy.  By Corollary~\ref{extend}, such foliations can be smoothly extended to a disk foliation on $Q_{\gamma}(\tau,A)$, and hence $V_{\gamma}(\tau,A)$, using the original $C^{\infty,0}$ foliation as a guide.

To simplify the exposition, we will fix a Riemannian metric on $M$ as described in Notation~\ref{metricmetric}.  

Recall that if two curves intersect in $\partial_v X$, for some codimension zero submanifold $X$ with piecewise horizontal and vertical boundary, the curve with greater slope, when viewed from outside of $X$, is said to {\sl dominate} the other curve.

\begin{notation} We say one object is $O(\epsilon)$ close to another if for some constant $K$ independent of the two objects, the objects are $K\epsilon$ close.
\end{notation}

\begin{thm} \label{C0V} Let $\mathcal F_0$ be a $C^{\infty,0}$ foliation of $M$, and let $(P,\mathcal P)$ and $g=g(P)$ be given as in Notation~\ref{metricmetric}.  Fix $\epsilon>0$ and let $V$, $N_h$, and $\mathcal F$ be given as in Theorem~\ref{metriccreateholonbhds}.  

 Then there are a regular neighborhood $N_v \subset V$ of the vertical edges of $\partial Q$ in $V$ and smooth plane fields $\xi_V^{\pm}$ defined on $V$ satisfying
\begin{enumerate}

\item $\xi_V^+$ is positive and $\xi_V^-$ is negative,

\item $\xi_V^{\pm}=T\mathcal F $ on $\overline{N_h\cup N_v}$ and is contact at all other points of $V$,

\item $\xi_V^+$ dominates $\mathcal F$ along $\partial_v V$, with the domination strict outside $\overline{N_h\cup N_v}$,

\item $\xi_V^-$ is dominated by $\mathcal F $ along $\partial_v V$, with the domination strict outside $\overline{N_h\cup N_v}$,

\item each of $\xi_V^{\pm}$ is positively transverse to $\Phi$, and

\item each of $\xi_V^{\pm}$ is $O(\epsilon)$ $C^0$ close to $T\mathcal F$ on $V$.

\end{enumerate}
\end{thm}

The proof of this theorem will occupy the rest of this section.  By symmetry, it will suffice to establish the existence of $\xi^+$.

It suffices to consider the case that $n=1$; so, to simplify notation, write $V=V_{\gamma}(\tau,A)$, with metric product neighborhood $(P,\mathcal P)$.  Recall that the metric product neighborhood $(P,\mathcal P)$ has product metric and horizontal product foliation $\mathcal P$ induced by the identification 
$P\cong [-1,1]\times S^1\times [-1,1]$.  Use this identification to view $V=V_{\gamma}(\tau,A)\subset P$ as a subset $V\subset [-1,1]\times S^1\times [-1,1]$.  Notice that $\partial_v V\subset \partial_v P$.  Let $N_R$ denote the open neighborhood of $R_{\gamma}(\tau,A)$ in $V$ given by the intersection of $V$ with $[-1,1]\times ([1/2,1]\cup [-1,-1/2])\times [-1,1]$.  Recall that $\mathcal F=\mathcal P$ on $N_R$.

Eventually, we will further constrain $N_v$, but for now, let $N_v\subset N_R$ be any regular neighborhood of the vertical edges of $\partial Q$, and set $N= N_h\cup N_v$.

Write $Q=Q_{\gamma}(\tau,A)$, and let $$\pi: Q\to V$$ denote the quotient map which reverses the splitting of $V$ along $R_{\gamma}(\tau,A)$.  Viewing $S^1$ as the quotient $[0,1]/0\sim 1$, the identification $V\subset [-1,1]\times S^1\times [-1,1]$ induces an identification
$$Q\subset [-1,1]\times [0,1]\times [-1,1]\subset \mathbb R^3,$$ with $\partial_v Q\subset \partial_v([-1,1]\times [0,1]\times [-1,1])$.  We will abuse notation and let $N_h$, $N_v$, $N$, and $N_R$ also denote their pullbacks to $Q$ under $\pi:Q\to V$.  Similarly, we let $\mathcal F$ denote the pullback $\pi^{-1}(\mathcal F\cap V)$ when this meaning is clear from context.

Our goal is to construct a smooth positive confoliation on $V$ satisfying the conditions of Theorem~\ref{C0V}.  We will do this by defining a smooth positive confoliation $\xi^+$ on $Q$ which smoothly glues, via $\pi: Q\to V$, to a smooth confoliation on $V$.  As a first step, we will define a smooth foliation on $\partial_v Q$ which will serve as the characteristic foliation of the contact structure $\xi^+$.  This characteristic foliation will be closely related to the restriction of $\mathcal F$ to $\partial_v Q$.  The following proposition will be used to make the transition from continuous to smooth structures.  

\begin{prop}\label{approxlem} Let $\Phi$, $N_h$, $N_v$, $N$, and $Q$ be as given above, and let $\kappa\in\{\pm 1\}$.

Let $X$ denote either a vertical face or a union of three vertical faces of $Q$, with $\partial_v X=\sigma\cup\tau$.  Let $\mathcal G_0$ be a $ C^{\infty,0}$ foliation of $X$ which is everywhere transverse to $\Phi$, satisfies $\mathcal G_0=\mathcal P$ on $\overline{N_v}\cap X$, and is smooth when restricted to $\overline{N}\cap X$.  Let $G_0:\sigma\to\tau$ be the holonomy map defined by following leaves of $\mathcal G_0$ across $X$, beginning in $\sigma$ and ending in $\tau$.  Let $\varepsilon:\sigma\to\mathbb R$ be a continuous function satisfying $\varepsilon(z)\ge 0$, with equality if and only if $z\in\partial \sigma$.

Then there is a $C^\infty$ foliation $\mathcal G$ on $X$, with holonomy across $X$ given by the holonomy map $G:\sigma\to \tau$, such that
\begin{enumerate}

\item $\mathcal G$ is positively transverse to $\Phi$,
\item $\mathcal G$ is arbitrarily $C^0$ close to $\mathcal G_0$,
\item $\mathcal G=\mathcal G_0$ on $\overline{N}$, 
\item $T\mathcal G$ dominates (respectively, is dominated by) $T \mathcal G_0$, with the domination strict outside $\overline{N}$, if $\kappa=1$ (respectively, $\kappa=-1$), and,
\item $|G(z)-G_0(z)|\le \varepsilon(z)$.
\end{enumerate}
\end{prop}
 
\begin{proof}
It suffices to consider the case that $\kappa=1$.  Begin by considering the case that $X$ lies in the vertical plane $x=1$.  

At each point $(1,y,z)$ of $X$ let $f(y,z)$ be the continuous function such that the line field $\partial_y + f(y,z) \partial_z$ is tangent to $\mathcal G_0$.  Let $g: X\to \mathbb R$ be a smooth function arbitrarily close to $f$ and such that $g(y,z)\ge f(y,z)$, with equality if and only if $(1,y,z)\in\overline{N}$.  Let $\mathcal G$ be the foliation given by the integral curves of the flow tangent to $\partial_y + g(y,z) \partial_z$.  Denote the foliation determined by $g$ by $\mathcal G$.  Certainly, conclusions (1)--(4), are satisfied by $\mathcal G$.

Let $\varepsilon_0$ be one half of the minimum value of $\varepsilon$ on $\sigma - \text{int}(N(\partial_hV))$, where $N(\partial_hV)$ is the portion of $N$ corresponding to a regular neighborhood of $\partial_hV$.

Let $g_n:X\to\mathbb R$ be a sequence of smooth functions, each determining a foliation satisfying (1)--(4) and with limit $f$.  By the smoothness of solutions to ODEs, there is a number $m$ so that by setting $g=g_n$ for any $n> m$, and letting $\mathcal G$ be the foliation determined by $g$, the corresponding holonomy map $G:\sigma\to\tau$ satisfies $G(z)-G_0(z)< \epsilon_0$.  Thus (5) holds for $z \notin \sigma - \text{int}(N(\partial_hV))$.  On $\text{int}(N(\partial_hV))$, $g=f$ , and hence $G(z)=G_0(z)$ for all $z\in \sigma \cap \text{int}(N(\partial_hV))$.  Thus (5) holds for all $z\in \sigma$.  

Finally, we consider the remaining possibilities for $X$.  Certainly the proof as given applies to any single face of $Q$.  In the case of a union of three faces of $Q$, isometrically flatten out the union so that it lies in a single plane to see that the proof as given applies.
\end{proof}

\begin{figure}[htbp] 
\centering
\includegraphics[width=5in]{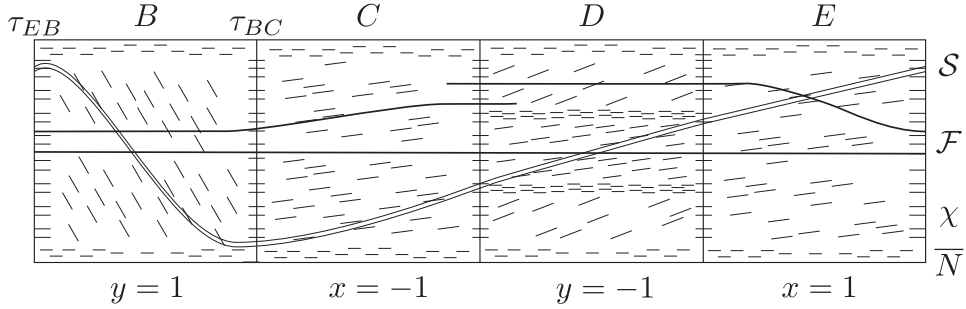} 
\caption{The vector field $\chi$ is shown as short dashes.  Two leaves of $\mathcal F$ are shown, one of which corresponds to the annulus leaf $A$.  A pair of circle leaves of $\mathcal S$ is shown dominating $\chi$.}
\label{thecube}
\end{figure}

Next we introduce some useful notation.  Label the vertical faces of $Q$ by $B,C,D,E$, where $\pi(B)\subset \pi(D)=R_{\gamma}(\tau,A)$, and the sequence $B,C,D,E$ is a listing of the faces in counterclockwise order about $\partial_v Q$.  These faces are illustrated in Figure~\ref{thecube}.  

Projecting these labels to $V$, we will abuse notation when convenient by considering $B$ to be a subset of $D$.  Notice that each of $B\cap \overline{N_h}$, $C\cap \overline{N_h}$, and $E\cap \overline{N_h}$ consists of two components, whereas $D\cap \overline{N_h}$ consists of four components.  

Also, label the vertical edges of $Q$:
$$\tau_{BC}=B\cap C,\, \tau_{CD}=C\cap D, \,\tau_{DE}=D\cap E,\,\mbox{and } \tau_{EB}=E\cap B.$$
Again, projecting these labels to $V$, we will abuse notation when convenient by considering $\tau_{EB}$ to be a subset of $\tau_{DE}$ and $\tau_{BC}$ to be a subset of $\tau_{CD}$.

Now consider the (continuous) holonomy maps of $\mathcal F$ across the 2-dimensional faces of $\partial_v Q$:
$$\tau_{EB}\overset{F_B}\longrightarrow \tau_{BC} \overset{F_C}\longrightarrow \tau_{CD}\overset{F_D}\longrightarrow \tau_{DE} \overset{F_E}\longrightarrow \tau_{EB}.$$
Since $\mathcal F$ is smooth on $N_h$, the restrictions of $F_B,F_C,F_D,$ and $F_E$ to $\tau_{BC}\cap N_h,\tau_{CD}\cap N_h,\tau_{DE}\cap N_h$ and $\tau_{EB}\cap N_h$, respectively, are smooth functions.

Recall that $\mathcal F=\mathcal P$ along $R_{\gamma}(\tau,A)$.  Therefore, using the identifications $V\subset [-1,1]\times S^1\times [-1,1]$ and $Q\subset [-1,1]\times [0,1]\times [-1,1]\subset \mathbb R^3$, $B\cong (A\cap B)\times \tau_{EB}$, with the leaves of $\mathcal F\cap B$ identified with the leaves $(A\cap B)\times \{t\}$, and $D\cong (A\cap D)\times \tau_{CD}$ with the leaves of $\mathcal F\cap D$ identified with the leaves $(A\cap D)\times \{t\}$.  In particular, $\tau_{EB}\cong\tau_{BC}$ and $\tau_{CD}\cong\tau_{DE}$, and under these smooth identifications, the maps $F_B$ and $F_D$ are identity maps.

In the next corollary, we show that the $C^{\infty,0}$ restrictions of $\mathcal F$ to $C$ and $E$ can be approximated by smooth foliations which dominate.

\begin{cor} \label{wigglysides}
There are smooth foliations $\chi_C$ on $C$ and $\chi_E$ on $E$ satisfying:
\begin{enumerate}

\item each of $\chi_C$ and $\chi_E$ is positively transverse to $\Phi$,
\item $\chi_C$ (respectively, $\chi_E$) is arbitrarily $C^0$ close to $\mathcal F$ on $C$ (respectively, on $E$),
\item $\chi_C = \mathcal F$ and $\chi_E = \mathcal F$ on $\overline{N}$,
\item $\chi_C$ (respectively, $\chi_E$) dominates $\mathcal F$ on $C$ (respectively, $E$), with the domination strict outside $\overline{N}$, 
\end{enumerate}
\end{cor}

\begin{proof} Let $\mathcal G_0$ denote the restriction of $\mathcal F$ to $X$, where $X$ is either $C$ or $E$, and set $\kappa=1$.  Apply Proposition~\ref{approxlem} with $X=C$ to obtain $\chi_C$, and with $X=E$ to obtain $\chi_E$.  
\end{proof}

In the next proposition, we show that the $C^{\infty,0}$ restrictions of $\mathcal F$ to $B$ and $D$ can be approximated by smooth foliations which dominate along $D$ and are dominated along $B$.  Moreover, the union of these foliations with the smooth foliations of the preceding corollary gives a foliation realizing decreasing monodromy about $\partial_v Q$.  

\begin{prop}\label{decrexists} Let $\chi_C$, $\chi_E$ be as guaranteed in Corollary~\ref{wigglysides}, with monodromy maps
$K_C:\tau_{BC}\to\tau_{CD}$ and $K_E:\tau_{DE}\to\tau_{EB}$ respectively.
There are smooth foliations $\chi_B$ and $\chi_D$, defining holonomy maps $K_B:\tau_{EB}\to\tau_{BC}$ and $K_D:\tau_{CD}\to\tau_{DE}$ respectively, satisfying 
\begin{enumerate}
\item when restricted to $\overline{N}$, $\chi_B=\mathcal F$ and $\chi_D=\mathcal F$,
\item $\mathcal F $ is dominated by $\chi_B$ along $B$, with the domination strict outside $\overline{N}$,
\item $ \mathcal F$ dominates $\chi_D$, with the domination strict outside $\overline{N}$,
\item $\chi_B$ and $\chi_D$ agree where identified by $\pi:Q\to V$,
\item $K_D(t)=K_B^{-1}(t)$ for all $t\in \tau_{EB}\subset\tau_{DE}$,
\item when restricted to $\overline{N}$, $K_B=F_B$ and $K_D=F_D$,
\item $K_B(t)\le F_B(t)$ for all $t\in \tau_{EB}$,
\item $K_D(t)\ge F_D(t)$ for all $t\in \tau_{CD}$,
\item $K_E K_D K_C K_B(t)\le t$ for all $t\in \tau_{EB}$, 
\item the foliation on $\partial_v Q$ defined by union
$$\chi = \chi_B\cup\chi_C\cup\chi_D\cup \chi_E$$
is smooth, and 
\item each of $T\mathcal K_B$ and $T\mathcal K_D$ is $O(\epsilon)$ $C^0$ close to $T\mathcal F$.  

\end{enumerate}
where, for each inequality, equality holds if and only if $t\in \overline{N}$.
\end{prop}

\begin{proof} 
Recall that the metric product neighborhood $(P,\mathcal P)$ has product metric and horizontal product foliation $\mathcal P$ induced by the identification $P\cong[-1,1]\times S^1\times [-1,1]$.  Use this identification to view $V_{\gamma}(\tau,A)\subset P$ as a subset $V_{\gamma}(\tau,A)\subset [-1,1]\times S^1\times [-1,1]$ and therefore $Q$ as a subset of $[-1,1]\times [0,1]\times [-1,1]$.  This identification induces $\tau_{CD}\cong\tau_{DE}\cong [-a,d]$ and $\tau_{EB}\cong\tau_{BC}\cong [-b,c]$ for some $-1\le -a<-b<0<c<d\le 1$.  Choose $e$ so that $c<e<d$ and $(c,e]$ is disjoint from $\overline{N}$.

First define $K_B$ and $K_D$ satisfying (5)--(9).  We do so as follows.  As a first step towards choosing $K_D$, let $\hat{K}_D:[-a,d]\to [-a,d]$ be any orientation preserving diffeomorphism that maps $[-b,c]\to [c,e]$ and satisfies $F_D(t) \le \hat{K_D}(t)$ for all $t\in \tau_{CD}$, with equality if and only if $t\in\overline{N}$.

Next choose $K_B$ satisfying (6), $K_B(t) \le (K_E \hat{K}_D K_C )^{-1}(t)$, and (7).  Finally, choose $K_D$ satisfying (6) to agree with $K_B^{-1}$ on $[-b,c]$ and so that $F_D(t) \le K_D(t)\le \hat{K}_D(t)$ for all $t\in [-a,-b]\cup[c,d]$.  Thus, $K_D(t)\le \hat{K}_D(t)$ for all $t\in [-a,d]$, and $K_D$ satisfies (8).  Since $K_B(t)<(K_E \hat{K}_D K_C )^{-1}(t) \le (K_E K_D K_C )^{-1}(t)$, condition (9) is satisfied.  Some of these relationships are shown in Figure~\ref{Monodromy}.

\begin{figure}[htbp] 
\centering
\includegraphics[width=4.5in]{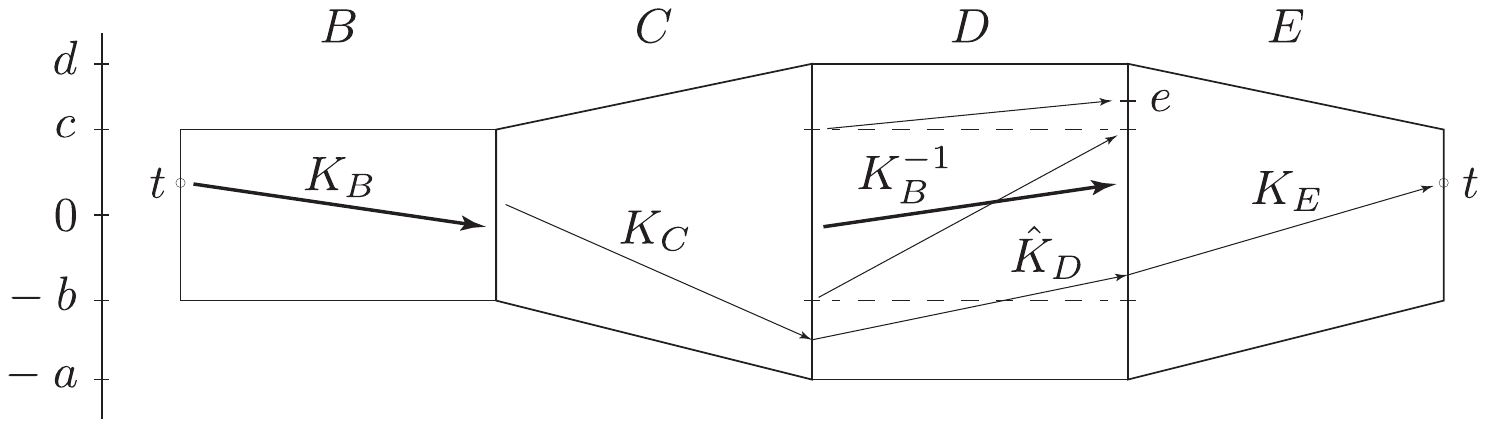} 
\caption{}
\label{Monodromy}
\end{figure}

Next, use $K_B$ and $K_D$ to construct $\chi_B$ and $\chi_D$ respectively.  By construction, $\mathcal F=\mathcal P$ in a neighborhood of $R_{\gamma}(\tau,A)$, and hence $\mathcal F$ is equal to $\mathcal P$ in a neighborhood of $B$ and $D$ in $Q$.  So $\mathcal F$ restricts to foliations of $B=[-1,1]\times \{1\}\times [-a,d]$ and $D=[-1,1]\times \{-1\}\times [-a,d]$, respectively, by horizontal straight line segments.  

The foliation of $B$ by line segments with endpoints $(1,1,z)$ and $(-1,1,K_B(z))$ has leaves given by $(1-\rho)(1,1,z) + \rho(-1,1, K_B(z))$ with $\rho \in [0,1]$ and has the desired monodromy $K_B(z)$.  To guarantee (1) holds, replace $\rho$ by $\rho(t)$ where $\rho:[-1,1] \to [0,1]$ is a smooth damping function chosen so that $\rho^{-1}(0) \cup \rho^{-1}(1)$ corresponds to $\overline{N} \cap B$.  Let $\chi_B$ be this damped linear foliation.  Similarly define $\chi_D$ to be the damped linear function with holonomy $K_D(z)$ damped so that leaves are horizontal exactly on $\overline{N} \cap D$.

The foliations $\chi_B$ and $\chi_D$ are smooth since the corresponding holonomy maps $K_B$ and $K_D$ are smooth.  In addition, since the foliations $\chi_B, \chi_C,\chi_D$ and $\chi_E$ are horizontal in a neighborhood of $\partial_v Q^{(1)}$, they glue together to give a smooth foliation of $\partial_v Q$.  Moreover, since $V$ is $\epsilon$-flat, the straight lines used to create $\chi_D$ have slope between $\pm \epsilon$, thus after damping, $T\chi_D$ has slope bounded in absolute value by $2\epsilon$.  Since $\mathcal F$ was chosen to satisfy the conclusions of Theorem~\ref{metriccreateholonbhds}, in particular, $\mathcal F$ is $\epsilon$-horizontal, it follows that $T\chi_D$ and $T\chi_B$ are $O(\epsilon)$ close to $T\mathcal F$.  
\end{proof}

Notice that $\chi$ is the characteristic foliation of a smooth 2-plane field along $\partial_v Q$ defined as follows.  At each point $p$ of $\partial_v Q$, let $\xi_p$ denote the 2-plane perpendicular to $\partial_v Q$ which contains $T_p\chi$.  We will show that this 2-plane field extends to a smooth confoliation on $V$ which stays close to $\mathcal F$.  The first step in constructing this extension is to build a circle foliation dominated by $\chi$ which in turn bounds a disk foliation of $Q$.

\begin{cor}\label{circfol} Let $K_B, K_C, K_D, K_E,$ and $\chi$ be as given in Proposition~\ref{decrexists}.  There is a smooth foliation $\mathcal S$ of $\partial_vQ$ by circles (with corners along $\partial_v Q^{(1)}$) such that 
\begin{enumerate}
\item $\chi$ dominates $T\mathcal S$, with the domination strict outside $\overline{N_h}$, 

\item $\mathcal S=\chi=\mathcal F$ on $\overline{N_h}$, and 

\item $\mathcal S$ is $C^0$ $O(\epsilon)$ close to each of $\mathcal F$ and $\chi$.  

\end{enumerate}
\end{cor}

\begin{proof} Begin by constructing $\mathcal S$ on $X$, where $X$ is the union of the faces $C, D,$ and $E$ of $Q$.  Let $\mathcal G_0$ denote the restriction of $\chi$ to $X$, and let $G_0$ denote the holonomy map $G_0:\tau_{BC}\to\tau_{EB}$ given by the composition $G_0= K_EK_DK_C$.  For each $z\in\tau_{BC}=[-b,c]$, let $\varepsilon(z)=K_B^{-1}z-G_0(z)$.  Apply Proposition~\ref{approxlem} to get a foliation $\mathcal G$ on $X$ satisfying
\begin{enumerate}

\item $\mathcal G$ is positively transverse to $\Phi$,
\item $\mathcal G$ is $\epsilon$ $C^0$ close to $\mathcal G_0$,
\item $\mathcal G=\mathcal G_0$ on $\overline{N}$, 
\item $T\mathcal G$ dominates $T \mathcal G_0$, with the domination strict outside $\overline{N}$ and,
\item $ G(z)-G_0(z)\le \varepsilon(z)$.  
\end{enumerate}
Since $\mathcal G_0$ is $O(\epsilon)$ close to $\mathcal F$ by Corollary~\ref{wigglysides} and Proposition~\ref{decrexists}, and $\mathcal G$ is $\epsilon$ close to $\mathcal G_0$, $\mathcal G$ is $O(\epsilon)$ close to $\mathcal F$.  Set $\mathcal S=\mathcal G$ on $X$.

Next, construct $\mathcal S$ on the remaining vertical face, $B$.  Take advantage of the fact that, as described in the proof of Proposition~\ref{decrexists}, $\chi$ restricted to $B$ consists of damped straight line segments.  Let $\mathcal S$ consist of a smooth family of similarly damped line segments from $G^{-1}(z)\in \tau_{EB}$ to $z\in\tau_{BC}$.  Condition (5) guarantees that $G(z)\le K_B^{-1}(z)$, with equality only in the collar of $\{ -b,c\}$ in $[-b,c]$ determined by $\overline{N}$, and hence line segments from $G(z)$ to $z$ are steeper than line segments $K_B^{-1}(z)$ to $z$ for $z$ outside this collar of $\{ -b,c\}$.  This property is preserved under suitably chosen damping.  Notice that $\mathcal S$ lies within $\epsilon$ of $\chi$ in $B$.

The foliation $\mathcal S$ is a foliation by circles since it has trivial monodromy about $\partial_v Q$.  By construction, $\mathcal S$ is dominated by $\chi$, with the domination strict exactly outside $\overline{N}$.  Moreover, $\mathcal S$ lies within $\epsilon$ of $\chi$, and hence $O(\epsilon)$ close to $\mathcal F$.
\end{proof}

Next we consider a smooth cylinder $Q'$ which lies nicely in $Q$.  We do this as follows.

Recall that in our preferred coordinates (namely, the ones inherited from $P$), $Q$ is an $\epsilon$-flat, $x$-invariant subset of $[-1,1]^3$ which is diffeomorphic to a cube and satisfies $[-1,1]\times [-1,1]\times \{0\}\subset Q$.  In particular, $\mathcal F$ agrees with the horizontal foliation $\mathcal P$ at all points $(x,y,z)$ of $Q$ for which $ y \notin [-1/2,1/2]$.

\begin{figure}[htbp] 
\centering
\includegraphics[width=5in]{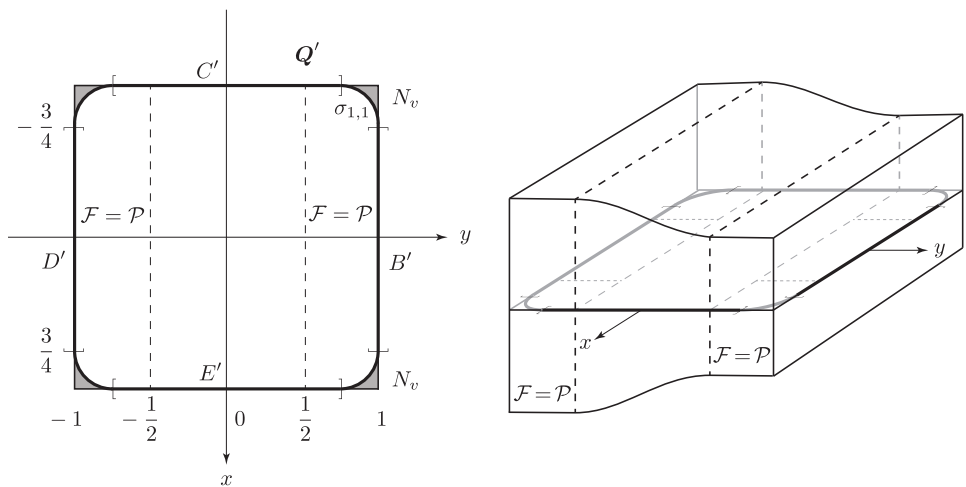} 
\caption{$\Delta\subset [-1,1]^2$ and $Q$.}
\label{Qprime}
\end{figure}

\begin{figure}[htbp] 
\centering
\includegraphics[width=4in]{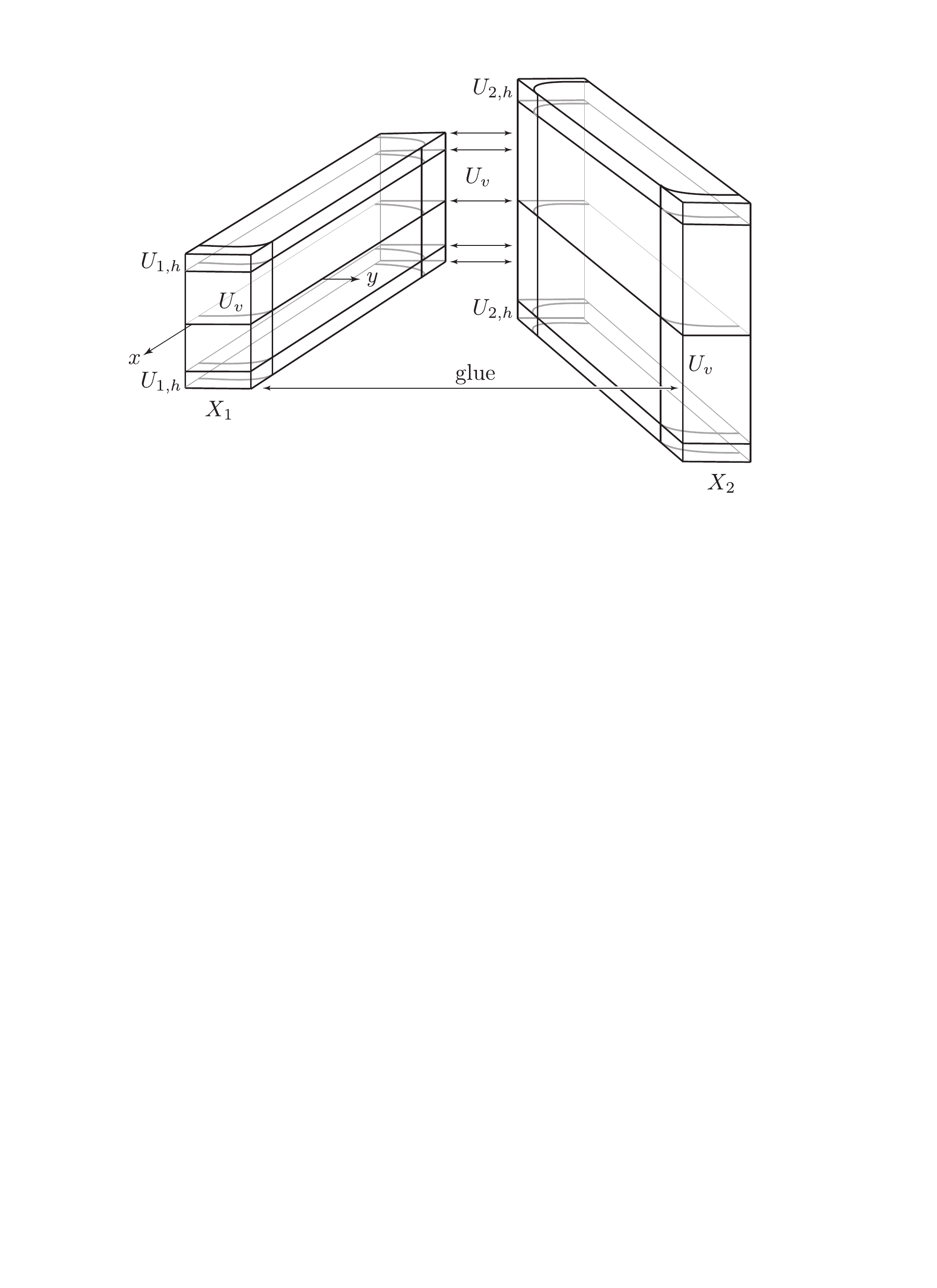} 
\caption{The glueing $X=X_1\cup X_2$.}
\label{gluingpic}
\end{figure}

Let $\Delta$ be the smooth disk embedded in the square $[-1,1]^2$ with smooth boundary 
$$\partial\Delta = B'\cup\sigma_{-1,1}\cup C'\cup\sigma_{-1,-1}\cup D'\cup\sigma_{1,-1}\cup E'\cup\sigma_{1,1},$$ where
$B' =[-3/4,3/4]\times \{1\}$,
$C'=\{-1\}\times [-3/4,3/4]$,
$D'=[-3/4,3/4]\times \{-1\}$, 
$E' =\{1\}\times [-3/4,3/4]$,
and $\sigma_{i,j}$ is a curve that rounds the corner near $(i,j)$ and smoothly connects the closest pair of line segments just defined.  See Figure~\ref{Qprime}.  Set $Q'=Q\cap (\Delta\times [-1,1])$, a smooth closed cylinder.  Finally, we specify a particular $N_v$: 
set $N_v=Q\setminus Q'$.  By the choice of $\Delta$, $N_v$ is an open regular neighborhood in $Q$ of the vertical edges of $Q$.

\begin{figure}[htbp] 
\centering
\includegraphics[width=2.4in]{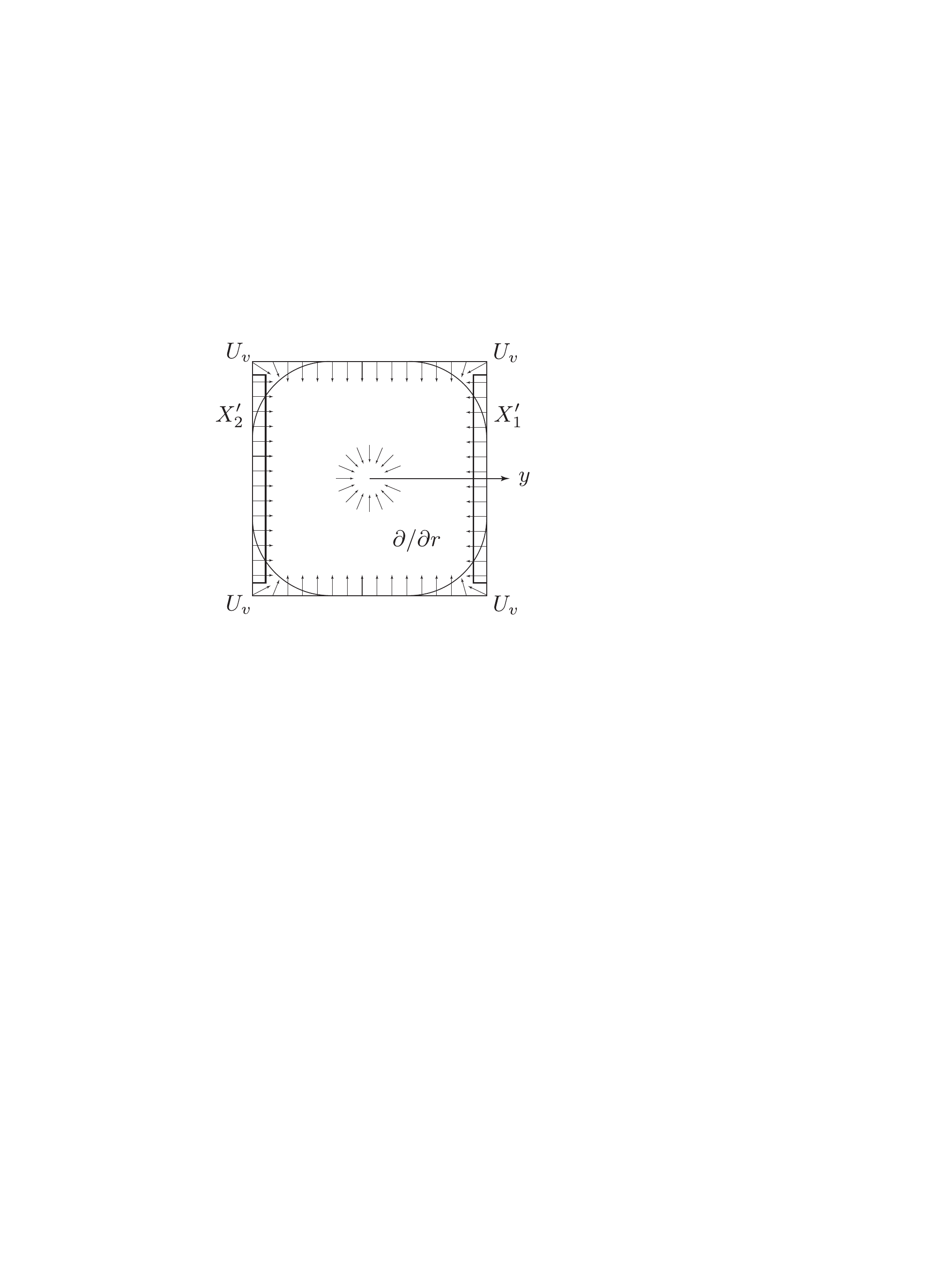} 
\caption{Radial flow on $\Delta$}
\label{radialflow}
\end{figure}

Next smoothly parametrize $\Delta$ by polar-like coordinates $(r, \theta)$, $r\in [0,1]$, where $\theta$ is the usual polar coordinate, and $r$ will be chosen to facilitate the identification of points $(x,-1)$ and $(x,1)$ in $\Delta$.  Let $X'_1=[-7/8,7/8] \times[7/8,1]$ and $X'_2=[-7/8,7/8] \times[-1,-7/8]$ be rectangular subsets of $\Delta$ as shown in Figure~\ref{radialflow}.  Choose $r$ so that 

\begin{enumerate}
\item when $(x,y)\in X'_1$, $r = -y$, 
\item when $(x,y)\in X'_2$, $r = y$, 
\item when $(x,y)\in \{1\}\times [-3/4,3/4]$, $\partial/\partial r=-\partial/\partial x$,
\item when $(x,y)\in \{-1\}\times [-3/4,3/4]$, $\partial/\partial r=\partial/\partial x$, and
\item the vector field $\partial/\partial r$ has a single, necessarily elliptic, singularity at $(0,0)$.
\end{enumerate}

Recall that $\chi$ is horizontal along $\partial_v Q\cap N_v$.  Let $\chi'$ denote the smooth extension of the line field $\chi$ on $\partial_v Q'\cap\partial_v Q$ to $\partial_v Q'$ obtained by defining $\chi$ on $\partial_v Q'\setminus\partial_v Q$ to be the tangent horizontal line field.  

Let $\mathcal S'$ denote a smooth extension of the foliation $\mathcal S\cap (\partial_v Q'\cap\partial_v Q)$ to a foliation by circles which satisfies
\begin{enumerate}
\item $\chi'$ dominates $T\mathcal S'$, with the domination strict outside $\overline{N_h}$, 

\item $\mathcal S'=\chi'=\mathcal F$ on $\overline{N_h}$, and 

\item $\mathcal S'$ is $C^0$ $O(\epsilon)$ close to each of $\mathcal F$ and $\chi'$.  

\end{enumerate} The existence of such an $\mathcal S'$ is guaranteed by continuity.

Next we extend the circle foliation $\mathcal S'$ to a smooth disk foliation of $Q'$.  

\begin{prop} \label{smooth disks} 
There is a smooth foliation $\mathcal D$ of $Q'$ by disks such that 
\begin{enumerate}
\item $T\mathcal D$ contains $\partial/\partial r$ in a width $1/4$ collar of $\partial_v Q'$,
\item $\mathcal D=\mathcal F$ on $\overline{N_h}$,
\item $\mathcal D$ is everywhere transverse to $\Phi$,
\item $\mathcal D$ is $O(\epsilon)$ $C^0$ close to $\mathcal F$, and 
\item $\mathcal S'=\mathcal D\cap \partial_v Q'$.
\end{enumerate}
\end{prop}

\begin{proof} 
Let $U$ denote the intersection of $Q'$ with the width $1/2$ collar of $\partial_v Q$ in $Q$.  In particular, $U$ contains the width $1/4$ collar of $\partial_v Q'$ in $Q'$.

Let $\mathcal H$ be a $C^{\infty}$ foliation of $Q'$ which is $\epsilon$ $C^0$ close to $\mathcal F$, is everywhere transverse to $\Phi$, is $x$-invariant, and satisfies $\mathcal H=\mathcal F$ on $N_h$ and $\mathcal H=\mathcal F=\mathcal P$ on $N_R\cap Q'$.

We will choose $\mathcal D$ to coincide with $\mathcal H$ outside $U$ and to smoothly interpolate between $\mathcal S'$ and $\mathcal H$ over $U$.  We will take advantage of the polar-like coordinates $(r,\theta)$ on $\Delta$.  Label the leaves of $\mathcal H$ by $H_t, t\in \tau$, where $H_t$ is the (disk) leaf intersecting $\tau$ at $t$, and let $h_t(r,\theta), (r,\theta)\in \Delta$, be the smooth family of functions such that the graph of $h_t$ is the leaf $H_t$.

Label the leaves of $\mathcal S'$ by $S'_t, t\in \tau$, where $S'_t$ is the circle leaf intersecting $\tau$ at $t$.  The leaf $S'_t$ can be described as a graph $z=s_t(\theta), \theta\in S^1$.  Since $\mathcal S'$ is smooth, $s_t$ defines a smooth family of smooth graphs.  Extend $\mathcal S'$ to a foliation $\tilde{\mathcal S'}$ of $U$, by extending the functions $s_t$ to functions $\tilde{s}_t$ defined on $U$ by $\tilde{s}_t(r,\theta) = s_t(\theta)$.  Since $\partial/\partial r$ lies in the restriction of $T\mathcal F$ to $U$, the leaves of $\tilde{\mathcal S'}$ lie in $Q'$ and describe a smooth foliation of $U$.  

Let $g$ denote a smooth bump function defined on $\Delta$ which is $1$ on a width $1/4$ collar of $\partial_v Q'$ and $0$ outside $U$.  Since $g$ is smooth and $U$ is compact, $g$ has bounded first partial derivatives, with the bounds independent of $\epsilon$.  Finally, define
$$d_t(r,\theta)=g(r,\theta)\tilde{s}_t(r,\theta)+(1-g(r,\theta)) h_t(r,\theta).$$ Let $\mathcal D$ be the smooth foliation with leaves given by the graphs of $z=d_t(r,\theta)$.  Computing first partial derivatives, we obtain: 
$$\frac{\partial d_t}{\partial r} = \frac{\partial g}{\partial r}\cdot (\tilde{s}_t-h_t) + (1-g)\frac{\partial h_t}{\partial r}$$ and 
$$\frac{\partial d_t}{\partial \theta} = \frac{\partial g}{\partial \theta}\cdot (\tilde{s}_t-h_t) + g\cdot (\frac{\partial \tilde{s}_t}{\partial \theta} - \frac{\partial h_t}{\partial \theta}) + \frac{\partial h_t}{\partial \theta}.$$

Since $Q'$ is $\epsilon$-flat, $|\tilde s_t - h_t| < \epsilon$.  The partials of $\tilde s_t$ and $\tilde h$ are $O(\epsilon)$ small since $\mathcal S', \mathcal H$, and $\mathcal F$ are $C^0$ $O(\epsilon)$ close to horizontal.  It follows that $\mathcal D$ is $O(\epsilon)$ close to horizontal, and hence $O(\epsilon)$ close to $\mathcal F$.  Finally, note that $\mathcal D=\tilde{\mathcal S}$ when $g=1$.  Hence, $\partial \mathcal D= \mathcal S'$, and $T\mathcal D$ contains $\partial/\partial r$ in a width $1/4$ collar of $\partial_v Q'$.
\end{proof}

We will use the foliation $\mathcal D$ to extend the line field $\chi'$ to a contact structure across $Q'\setminus N_h$, and thus to define a smooth confoliation on $Q$ which is contact on the complement of $\overline{N_h\cup N_v}$.  First, we establish an elementary glueing lemma.

\begin{lemma} \label{stdsmoothing}
Suppose $X$ decomposes as a union of two cubes $X_1$ and $X_2$, where for some $u<v<w$ and nondegenerate closed intervals $J_1\subset J_2$, $X_1 = [-1,1]\times [u,v]\times J_1$ and $X_2 = [-1,1]\times [v,w]\times J_2$.  Let $$\alpha_1= dz-a(x,y,z) dx $$ 
be a smooth 1-form defining a positive confoliation $\xi_1= ker \alpha_1$ on $X_1$, and let $$\alpha_2= dz-b(x,y,z) dx $$ be a smooth 1-form defining a positive confoliation on $X_2$.  Let $U_{i,h}$ be a regular open neighborhood of $\partial_h X_i$ in $X_i$, $i=1,2$, and let $U_v$ be a regular open neighborhood of the faces $x=\pm 1$ in $X$.  Suppose that $U_{2,h} \cap X_1 \subset U_{1,h}$ (this allows $U_{2,h} \cap X_1=\emptyset$).  

In addition, suppose that the functions $a$ and $b$ satisfy the following:

\begin{enumerate}

\item $a=b$ on $X_1 \cap X_2$,
\item $a(x,y,z)=0 \iff a_y(x,y,z) =0 \iff (x,y,z)\in \overline{U_v \cup U_{1,h}}$, and 
\item $b_y(x,y,z) =0 \iff (x,y,z)\in \overline{U_v \cup U_{2,h}}$.  
\end{enumerate}

Then there is a smooth 1-form $\alpha = dz-c(x,y,z) dx$ defining a positive confoliation $\xi=\mbox{ ker }\alpha$ on $X$, where $c$ satisfies

\begin{enumerate}
\item $c=a$ on a neighborhood of the $y=u$ face of $X_1$,
\item $c=b$ on a neighborhood of the $y=w$ face of $X_2$, and
\item $c(x,y,z) = 0 \iff c_y(x,y,z) = 0 \iff (x,y,z)$ is in the closure of $U_v \cup U_{1,h} \cup U_{2,h}$.  
\end{enumerate}

Moreover, if each of $a$ and $b$ is $C^0$ close to $0$, then so is $c$.  
\end{lemma}

In other words, the continuous 1-form $\alpha_1\cup \alpha_2$ on $X$ can be $C^0$ approximated by a smooth 1-form $\alpha$ which agrees with $\alpha_1\cup \alpha_2$ on a neighborhood of $\partial X$, and describes a positive confoliation $\xi= \text{ker}(\alpha)$ on $X$ which is a contact structure exactly where $\xi_1$ or $\xi_2$ is a contact structure.

\begin{proof} For $i=1$ or $2$, $\xi_i$ is a positive confoliation, that is, $\alpha_i\wedge d\alpha_i\ge 0$.  Hence $a_y\ge 0$ and $b_y\ge 0$.  Also, by hypothesis, $a=b$ on $X_1 \cap X_2$.  Therefore, for each $(x_0,z_0)\in [-1,1]\times J_1$, the one variable functions $a(x_0,y,z_0), y\in [u,v]$, and $b(x_0,y,z_0), y\in [v,w],$ piece together to give a continuous function defined on $[v,w]$ which is smooth on the complement of $\{v\}$.  

Notice that $b(x,y,z)=0$ if and only if $(x,y,z)\in U_{1,h}\cap X_2$ or $b_y(x,y,z)=0$.  So $b(x,y,z)=0$ if and only if $(x,y,z)\in (U_{1,h}\cap X_2)\cup \overline{U_v \cup U_{2,h}}$.

To facilitate blending $a$ and $b$ into a smooth function $c$, choose a smooth function $\tilde b$ on $X_2$ such that: 
\begin{enumerate}
\item $\tilde b = b$ in a neighborhood of the $y=w$ face of $X_2$, 
\item $\tilde b_y \ge 0$
\item $\tilde b \ge b$ on $X_2$
\item $\tilde{b}(x,v,z)>b(x,v,z)\iff (x,v,z)\notin \overline{U_v\cup U_{2,h}}$
\item $\tilde{b}=0\iff$$ \tilde b_y=0 \iff b_y=0$ 
\end{enumerate}

At a point $(x,v,z)\in X_1 \cap X_2$ it might be that $a_y(x,v,z) > b_y(x,v,z) >0$, but in this case, the choice of $\tilde b$ forces $a(x,v,z) < \tilde b(x,v,z)$.  Thus it is possible to define a smooth $\tilde a$ on $X_1 \cup X_2$ which satisfies all of the following, and in particular (3): 
\begin{enumerate}
\item $\tilde{a}=a$ on $X_1$,
\item $\tilde a_y \ge 0$,
\item $\tilde a \le \tilde b$, and 
\item $\tilde a(x,y,z)=0 \iff \tilde a_y(x,y,z) =0 \iff (x,y,z)$ is an element of the closure of $U_v \cup U_{1,h} \cup U_{2,h}$.  
\end{enumerate}

To produce such an $\tilde a$, pick a non-negative smooth extension $e(x,y,z)$ of $a_y(x,y,z)$ to $X_1 \cup X_2$ which is $0$ exactly on $\overline{U_v \cup U_{1,h} \cup U_{2,h}}$.  With (3) in mind, such an extension can be modified by multiplying by a smooth non-negative map which takes the value $1$ on $X_1$ and approaches $0$ quickly on $X_2$ so that 
$$\tilde a(x,y,z) = \int_u^y e(x,s,z)ds + a(x,u,z)$$
satisfies (1)--(4).

Now choose $\sigma:[u,w] \to [0,1]$ such that $\sigma([u,v]) = 0$, $\sigma([\frac{v+w}{2},w]) = 1$, and $\sigma$ maps $(v,\frac{v+w}{2} )$ diffeomorphically to $(0,1)$.  Set 
$$c(x,y,z)=(1-\sigma(y))\tilde a(x,y,z) + \sigma(y) \tilde b(x,y,z).$$
so that Properties (1) and (2) of $c$ are immediate.  Since
$$c_y=\sigma_y (\tilde b - \tilde a) + (1 - \sigma) \tilde a_y + \sigma \tilde b_y \ge 0,$$
is the sum of three non-negative terms, $\xi$ is a positive confoliation.  Moreover, if $c_y=0$ either $a_y=0$ and $(x,y,z)\in \overline{U_v \cup U_{1,h}}$, or $b_y=0$ and $(x,y,z)\in \overline{U_v \cup U_{2,h}}$.  Hence, if $c_y(x,y,z)=0$, necessarily $(x,y,z)\in \overline{U_v \cup U_{1,h} \cup U_{2,h}}$.  But this means also that $\tilde{a}=\tilde{b}=0$, and so $c=0$.

Conversely, suppose $(x,y,z) \in \overline{U_v \cup U_{1,h} \cup U_{2,h}}$.  It suffices to consider the case that $(x,y,z)\in U_v \cup U_{1,h} \cup U_{2,h}$, and hence at least one of the following is true:
\begin{enumerate}
\item $(x,y,z) \in U_v$ means $\tilde{a}=0$ and $\tilde{b}=b=0$,
\item $(x,y,z) \in U_{1,h}$ means $\sigma=0$ and $\tilde{a}=0$, and 
\item $(x,y,z) \in U_{2,h}$ means $\tilde{b}=b=0$ and $\tilde{a}=0$.  
\end{enumerate}

In each of these three cases, $c=0$ on an open set about $(x,y,z)$ and therefore $c_y(x,y,z)=0$.  Hence, property (3) of $c$ is satisfied.  
\end{proof}

\begin{cor}\label{extend} There exists a smooth confoliation $\xi$ on $Q$ that is $O(\epsilon)$ $C^0$ close to $T\mathcal F$, has characteristic foliation $\chi$ on $\partial_v Q$, and satisfies $\xi=T\mathcal F$ on $\overline{N_h\cup N_v}$.  Moreover, $\xi$ can be chosen so that $\pi(\xi)$ is a smooth confoliation on $V$, where $\pi$ is the quotient map $\pi:Q\to V$.
\end{cor}

\begin{proof} At all points of $\overline{N}=\overline{N_h\cup N_v}$, let $\xi$ be the tangent plane to $\mathcal F$.  Thus along $\overline{N} \cap \partial_v Q$, $\xi$ contains $\chi$, and along $\overline{N} \cap \partial_v Q'$, $\xi$ contains $\chi'$.

The foliation by disks, $\mathcal D$, given by Proposition~\ref{smooth disks} will be used to to extend $\xi$ to all of $Q'$.  

Let $\iota$ denote the smooth inward pointing vector field on $Q'$ given by lifting the vector field $-\partial/\partial r$ to the leaves of $\mathcal D$, where the lift is the pullback under the projection $(r,\theta,z)\to (r,\theta)$.  In particular, in the width $1/4$ collar about $\partial_v Q'$, $\iota=-\partial/\partial r$.

Notice that $\chi'$ and $\iota$ span a plane at every point of $\partial_v Q'$.  Denote this plane by $\xi$.  

At this point we have the start of a smooth confoliation $\xi$ on $\partial_v Q' \cup \overline{N}$, and $Q'$ is foliated by disks of $\mathcal D$ which are in turn either everywhere tangent to $\xi$ or foliated by a vector field $\iota$ which serves as a candidate for a Legendrian vector field.  This is directly analogous to a cylindrical neighborhood of the $z$-axis in the standard radial model of a tight contact structure.

To extend $\xi$ across $Q'$ it is enough to map $Q'$ to a standard model, use the technique of Lemma~5.14 of \cite{KR2}, and pull back the resulting contact structure to a contact structure $\xi$ on $Q'$.  Roughly speaking, $Q'$ is mapped to a solid cylinder in $\mathbb R^3$ centered along the $z$-axis in such a way that $\mathcal D$ and $\iota$ are mapped to horizontal planes and radial lines, respectively.  The standard radially symmetric contact structure is then pulled back to $Q'$.

Some extra care is needed so that the confoliation $\xi$ on $Q$ is $O(\epsilon)$ close to $T\mathcal F$.  By construction, it is $O(\epsilon)$ close along $\partial_v Q$ and {\it in coordinates} the planes of $\xi$ monotonically approach horizontal planes in $\mathbb R^3$ as you move radially towards the $z$-axis along Legendrian curves.  The issue is that the pullback confoliation planes may not monotonically approach $T \mathcal D$.

Since $Q'$ is compact, the metric distortion when compared to the standard model is bounded, and so it suffices to show that we can reduce to the case that the line field $\chi'$ is arbitrarily close to $T\mathcal S'$.  We do this by taking advantage of the width $1/2$ collar, $U$, of $\partial_v Q$ in $Q$ to define a contact structure with planes that rotate from slope $\chi'$ to slope close to $T\mathcal S'$ as follows.

Since the restriction of $\partial/\partial r$ to $U$ lies in both $T\mathcal F$ and $T\mathcal D$, the restrictions to $U$ of the flow lines for $\partial/\partial r$ lie in both $\mathcal D$ and $\mathcal F$.

 Recall that $\partial \Delta$ corresponds to $r=1$, and let $s(\theta,z)$ denote the slope of $T\mathcal S'$ at $(1,\theta,z)$.  Note that $s(\theta,z)$ is also the slope of $T\mathcal D\cap \partial_v U$ at $(r_0,\theta,z)$, where $r_0$ is determined by the condition that $(r_0,\theta,z)\in (\partial_v U\setminus \partial_v Q)$.
Let $c(\theta,z)$ denote the slope of $\chi'$ at $(1,\theta,z)$.  Since $\chi'$ strictly dominates $\mathcal S'$ exactly on the complement of $\overline{N_h}$, 
$$s(\theta,z)-c(\theta,z)\ge0, $$
with equality if and only if $(1,\theta,z)\in \overline{N_h}$.

Fix $\delta\in (0,1/2)$, and let $f: U \to [\delta,1]$ be a smooth function which satisfies 
\begin{enumerate}
\item $f(r,\theta,z)=1$ for all $(r,\theta,z)\in \overline{N}\cup \partial_v Q'$,
\item $f(r,\theta,z)=\delta$ for all $(r,\theta,z)\in (\partial_v U\setminus (\partial_v Q\cup N_h))$, and 
\item $f_r(r,\theta,z)\ge 0$ for all $(r,\theta,z)\in U$, with equality if and only if $(r,\theta,z)\in \overline{N}$.
\end{enumerate}

For $(r,\theta,z)\in U$, define 

$$ \alpha = dz-[(1 - f(r,\theta,z))s(\theta,z)+f(r,\theta,z)c(\theta,z)]d\theta.$$ 

 Setting $g(r,\theta,z)=(1 - f(r,\theta,z))s(\theta,z)+f(r,\theta,z)c(\theta,z)$, we have 

$$\alpha = dz - g(r,\theta,z) d\theta.$$

Notice that along $\partial_v Q'$, $\text{ker}(\alpha) = \xi$.  In addition, 
 $$d\alpha = -g_r dr d\theta - g_z dz d\theta,$$

and hence 
$$\alpha\wedge d\alpha = -g_r dr d\theta dz = -(f_r(r,\theta,z)(c(\theta,z)-s(\theta,z))) dr d\theta dz$$
is a smooth positive confoliation defined on $U$.  This confoliation agrees with $T\mathcal F$ on $U\cap \overline{N}$ and is a contact structure on $U\setminus \overline{N}$.  Moreover, by choosing $\delta$ as small as is necessary, we may guarantee that at each point $(r,\theta,z)\in (\partial_v U\setminus \partial_v Q)$ the line field given by $\xi$ restricted to $\partial_v U$
dominates and is as close to the line field given by $T\mathcal D$ restricted to $\partial_v U$ as is required.  Hence there is a smooth confoliation $\xi$ on $Q$ that is $O(\epsilon)$ $C^0$ close to $T\mathcal F$, has characteristic foliation $\chi$ on $\partial_v Q$, and satisfies $\xi=T\mathcal F$ on $\overline{N_h\cup N_v}$.  

To ensure that the contact structure on $Q$ glues smoothly to $V$, some further care is needed.  Extending the two plane field from the vertical boundary of a cylinder across a radially foliated disk involves a choice of rate of rotation of contact planes to horizontal.  The rates should be chosen to respect the glueing of $B$ to $D$.  Since $\xi=T\mathcal F$ along $N_v$, it suffices to show that the rates can be chosen to respect the glueing of $B'$ to $D'$.  To see that this is possible, proceed as follows.  

Recall the rectangles $X'_1, X'_2 \subset \Delta$, and for $i=1,2$, let $X_i$ denote the points $(r,\theta,z)\in Q$ with $(r,\theta)\in X_i'$.  Let $\xi_i$ denote the restriction of $\xi$ to $X_i$.  In terms of $(x,y,z)$ coordinates and given the form of $\xi$ in $U$, we can write $\xi_1= \text{ker}(dz-a(x,y,z) dx)$ and $\xi_2=\text{ker}(dz-b(x,y,z))$, for smooth functions $a(x,y,z)$ defined on $X_1=[-7/8,7/8]\times [7/8,1]\times \sigma$ and $b(x,y,z)$ defined on $X_2=[-7/8,7/8]\times [-1,-7/8]\times\tau$.

Without changing notation, regard $X_1$ and $X_2$ as subsets of $V$ using the quotient map $\pi:Q \to V$ given by identifying $y=\pm 1$.  Since $\xi_1=\xi_2$ along $X_1\cap X_2$, $a(x,y,z)=b(x,y,z)$ along $X_1\cap X_2$ in $V$.

Recall also that $\chi $ is dominated by $\mathcal F$ along $D$ and dominates $\mathcal F$ along $B$.  Hence, $a(x,y,z)=b(x,y,z)\ge 0$ along $X_1\cap X_2$, with equality if and only if $(x,y,z)\in \overline{N_h\cup N_v}$.  Lemma~\ref{stdsmoothing} therefore applies and produces a smooth contact structure on $X$ which agrees with $\xi_1\cup\xi_2$ in a neighborhood of $\partial X$.
\end{proof}

\begin{proof}[Proof of Theorem~\ref{C0V}] Let $\epsilon >0$.  By symmetry, it suffices to establish the existence of $\xi^+$.

Let $\mathcal F_0$ be a $C^{\infty,0}$ foliation of $M$, and let $\mathcal F$ be an $\epsilon$ $C^0$ close $C^{\infty,0}$ foliation constructed using Theorem~\ref{metriccreateholonbhds}.  Let $V(\tau, A)$ be one of the holonomy neighborhoods constructed for $\mathcal F$.  

Corollary~\ref{wigglysides} and Proposition~\ref{decrexists} guarantee the existence of a smooth foliation $\chi$ on $\partial_v Q$ which is dominated by and $O(\epsilon)$ $C^0$ close to $\mathcal F$.  Corollaries~\ref{circfol}--\ref{extend} show that $\chi$ is the characteristic foliation of a smooth confoliation on $Q$ that glues to a smooth confoliation on $V$.  This confoliation restricts to a contact structure outside $\overline{N}$.  Each of $\xi_V^{\pm}$ is $O(\epsilon)$ $C^0$ close to $\mathcal F$ on $V$, and hence $O(\epsilon)$ $C^0$ close to $T\mathcal F_0$ on $V$.  
\end{proof}

\section{Foliations with only trivial holonomy}\label{noholo}

In this section we investigate taut, transversely oriented $C^0$ foliations with only trivial holonomy.  First we recall some classical results, with a focus on smoothness assumptions.  New results then appear as Theorem~\ref{closetofibration}, Theorem~\ref{measuredimpliesfibered} and Corollary~\ref{allplanes}, with Corollary~\ref{allplanes} giving the conclusion needed for the main result in this paper.

We begin by recalling a theorem found in \cite{Im}.

\begin{thm} [Theorem~4.1 of \cite{Im}] \label{Imcover}
Let $\mathcal F$ be a taut, transversely oriented $C^0$ foliation in $M$ with only trivial holonomy, and let $\Phi$ be a $C^0$ flow transverse to $\mathcal F$.  Let $\tilde{M}$ denote the universal cover of $M$, and let $\tilde{\mathcal F}$ and $\tilde{\Phi}$ be the lifts of $\mathcal F$ and $\Phi$ to $\tilde{M}$.  Then $\tilde{M}$ is homeomorphic to $\mathbb R^2\times \mathbb R$, where each $\mathbb R^2\times \{z\}$, $z\in \mathbb R$, is a leaf of $\tilde{\mathcal F}$ and each $\{x\}\times \mathbb R$, $x\in\mathbb R^2$, is an orbit of $\tilde{\Phi}$.
\end{thm}

This theorem is given in \cite{Im} for foliations that are $C^{\infty}$.  However, this smoothness hypothesis in unnecessary.  For completeness, we include Imanishi's proof here, reframed using the language of leaf spaces and with careful attention paid to smoothness assumptions.

Begin by recalling a lemma found in \cite{Im}.  
\begin{lemma}\label{normalfence}\cite{Im}
Let $\mathcal F$ be a transversely oriented $C^0$ foliation with only trivial holonomy, and let $\Phi$ be a $C^0$ flow transverse to $\mathcal F$.  Suppose $H:[0,1]\times [0,1)\to M$ is a continuous map such that $H([0,1]\times \{t\})$ is a curve in a leaf of $\mathcal F$ for all $s$, and $H(\{s\}\times [0,1)$ is an immersed curve in a flow line of $\Phi$ for all $t$.  If $H$ extends continuously to $\{0\}\times [0,1]$, then $H$ extends continuously to $[0,1]\times [0,1]\to M$.
\end{lemma}

\begin{proof}
This follows immediately from Theorem~3.1 of \cite{Im}.  See also Lemma~9.2.4 of \cite{CC}.
\end{proof}

\begin{proof}[Proof of Theorem~\ref{Imcover}]
Let $\tilde{M}$ be the universal cover of $M$ and let $\tilde{\mathcal F}$ and $\tilde{\Phi}$ be the lifts of $\mathcal F$ and $\Phi$ to $\tilde{M}$.  Let $T$ be the leaf space of $\mathcal F$, and let $\rho:\tilde{M}\to T$ denote the associated quotient map.  Note that it is sufficient to prove that $\rho(\tilde{C})=T$ for any orbit $\tilde{C}$ of $\tilde{\Phi}$.

So suppose that $\rho(\tilde{C})\ne T$ for some orbit $\tilde{C}$ of $\tilde{\Phi}$.  Since $\tilde{C}$ is everywhere transverse to $\tilde{\Phi}$, $\rho(\tilde{C})$ is an embedded copy of $\mathbb R$ in $T$.  Note that although $\tilde{C}$ is properly embedded in $\tilde{M}$, $\rho(\tilde{C})$ may or may not be properly embedded in $T$.  In either case, there is a leaf $\tilde{L}$ of $\tilde{\mathcal F}$ such that
$\rho(\tilde{L})$ lies in the closure of $\rho(\tilde{C})$ but not in $\rho(\tilde{C})$.  Let $\tilde{C}'$ be an orbit of $\tilde{\Phi}$ passing through $\tilde{L}$.
\begin{figure}[htbp] 
\centering
\includegraphics[width=2.3in]{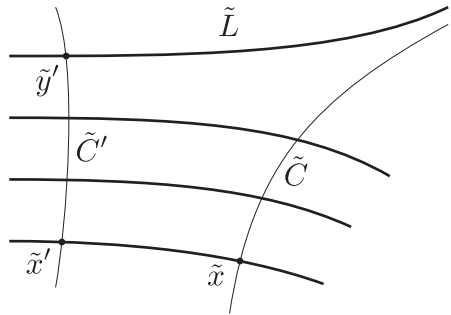} 
\caption{A continuous $H$ which does not continuously extend.  }
\label{Imanish}
\end{figure}
There is an interval $[x,y]$ in $T$ such that $[x,y]\cap \rho(\tilde{C}')=[x,y]$ and $[x,y]\cap \rho(\tilde{C})=[x,y)$.  As illustrated in Figure~\ref{Imanish}, 
this is impossible by Lemma~\ref{normalfence}.
\end{proof}

\begin{cor} \label{ImI}
Let $\mathcal F$ be a taut, transversely oriented $C^0$ foliation with only trivial holonomy, and let $\Phi$ be a $C^0$ flow transverse to $\mathcal F$.  Let $\Lambda$ be a minimal set of $\mathcal F$.  Then any region complementary to $\Lambda$ is a product $L\times [0,1]$, with each $\{x\}\times [0,1]$ a segment of the flow $\Phi$.  Moreover, the restriction of $\mathcal F$ to this complementary region is, up to a $\Phi$ preserving isotopy, the product foliation $L\times [0,1]$.
\end{cor}

\begin{proof}
Let $X$ denote the metric closure of some complementary region of $\Lambda$.  By Theorem~\ref{Imcover}, any lift of $X$ to $\tilde{M}$ has the form $\tilde{L}\times [0,1]$, where $\tilde{L}$ is the lift of a boundary leaf $L$ of $X$, and each $\{\tilde{x}\}\times [0,1]$ is a segment of the flow $\tilde{\Phi}$.  It follows that $X$ is an $I$-bundle, with each $I$-fiber a segment of a flow line of $\Phi$.  Now apply Lemma~\ref{noholIbundle} to conclude that the restriction of $\mathcal F$ to $X$ is, up to a $\Phi$ preserving isotopy, the product foliation $L\times [0,1]$.
\end{proof}

\begin{cor}\label{trivholo}
Let $\mathcal F$ be a taut, transversely oriented $C^{1,0}$ foliation with only trivial holonomy.  Then exactly one of the following is true:
\begin{enumerate}
\item $\mathcal F$ is a fibering of $M$ over $S^1$, 
\item $\mathcal F$ is minimal, or, 
\item $\mathcal F$ is a single Denjoy blow up of a minimal foliation; equivalently, $\mathcal F$ contains a unique minimal set, and this minimal set is exceptional with complement a product.
\end{enumerate}
\end{cor}

\begin{proof}
Let $\Phi$ be a $C^0$ flow transverse to $\mathcal F$.  If $\mathcal F$ is not minimal, then it contains a minimal set, $\Lambda$ say.  Let $X$ denote the metric closure of any component of the complement of $\Lambda$.  It follows from Corollary~\ref{ImI} that $X$ is an $I$-bundle $L\times [0,1]$, and the restriction of $\mathcal F$ to $X$ is, up to $\Phi$-preserving isotopy, the product foliation $L\times[0,1]$.  If $\Lambda$ is a compact leaf, conclude that $\mathcal F$ is a fibering of $M$.  Otherwise, $\Lambda$ is exceptional, and $\mathcal F$ is a single Denjoy blow up (along at most countably many leaves) of a minimal foliation; in other words, $\mathcal F$ contains a unique minimal set, and this minimal set is exceptional with complement a product.
\end{proof}

\begin{thm} [Theorem~1.3, \cite{Im}] \label{Imflow}
Let $\mathcal F$ be a taut, transversely oriented $C^0$ foliation with only trivial holonomy.  Suppose $\mathcal F$ does not contain an exceptional minimal set.  Then there is a topological flow $\psi:M\times \mathbb R\to M$ such that 
\begin{enumerate}
\item $\psi$ preserves $\mathcal F$; i.e., $\psi(-,t)$ sends each leaf of $\mathcal F$ into a leaf of $\mathcal F$.
\item $\psi$ is topologically transverse to $\mathcal F$; ie, $\psi(x,\mathbb R)$ is topologically transverse to $\mathcal F$.
\end{enumerate}
\end{thm}

Recall that a {\sl transverse measure} on a codimension one foliation $\mathcal G$ is an {\sl invariant} measure on each arc transverse to $\mathcal G$ that is equivalent to Lebesgue measure on an interval of $\mathbb R$.  Invariant, in this context, means that the measure of a transverse arc is unchanged under isotopies of the arc that keep each point on the same leaf of $\mathcal G$.

\begin{lemma}\label{measure} Let $\psi:M\times \mathbb R\to M$ be the topological flow of Theorem~\ref{Imflow}.  If $\psi(x_1,[s_1,t_1])$ and $\psi(x_2,[s_2,t_2])$ are isotopic through an isotopy that keeps each point on the same leaf of $\mathcal F$, then $t_1-s_1=t_2-s_2$.
\end{lemma}

\begin{proof} The isotopy can be lifted to the universal cover of $M$, so without changing notation we take $M=\tilde M$.  The advantage of working in $\tilde M$ is that leaves of the foliation are in bijective correspondence with $\psi(x_2,t)$ for $t\in \mathbb R$.

Let $\tau_i = \psi(x_i,[s_i,t_i])$.  The given isotopy sweeps out a family of curves $\alpha_t$ each contained in a leaf of $\mathcal F$ that starts at a point of $\tau_1$ and ends at a point of $\tau_2$, with $t\in [s_1,t_1]$ so that $\alpha_t(0)=\psi(x_1,t)$.  Another family of arcs $\beta_t$ each contained in a leaf of $\mathcal F$ can be generated by using the flow.  Define $\beta_t(s) = \psi(\alpha_{s_1}(s), t-s_1)$.  

It follows that $\beta_t(0)=\psi(\psi(x_1,s_1),t - s_1) = \psi(x_1, t) = \alpha_t(0)$ for $t\in [s_1,t_1]$, while $\beta_t(1) = \psi(\psi(x_2,s_2),t-s_1)= \psi(x_2, s_2 +t-s_1)$.

The arcs $\alpha_t$ and $\beta_t$ have the same initial point, are both contained in the same leaf, and terminate in the flow segment $\psi(x_2, [s_2+t_1 - s_1, t_2])$.  It follows that they terminate in the same leaf, that is, $s_2+t_1 - s_1 = t_2 $.
\end{proof}

\begin{cor}[Corollary~4.1, \cite{Im}]\label{measured} 
Let $\mathcal F$ be a taut, transversely oriented $C^{1,0}$ foliation with only trivial holonomy.  Suppose $\mathcal F$ does not contain an exceptional minimal set.  Then $\mathcal F$ admits a transverse measure.
\end{cor}

\begin{proof} Let $\tau$ be a transversal to $\mathcal F$.  First we show that if the terminal point of $\tau$ is $\psi(x,t)$, then there is an $[s,t] \subset \mathbb R$ such that $\tau$ is isotopic to $\psi(x,[s,t])$ through an isotopy that keeps each point on the same leaf of $\mathcal F$.

To see this, fix the initial point $x_0$ of $\tau$, and use the flow to homotop $\tau$ to an arc $\sigma$ starting at $x_0$ and contained in a leaf of $\mathcal F$.  This homotopy can be thought of as the image of a triangle $T$, with edges $\tau, \sigma$, and a flow arc $\psi(x,[s,t])$ for some value of $s$ swept out by $\psi(x,t)$.  The desired isotopy of $\tau$ is given by following horizontal arcs of $T \cap \mathcal F$.

Define the length of $\tau$ to be $|t-s|$.  Lemma~\ref{measure} guarantees that this defines a positive transverse measure on $\mathcal F$.
\end{proof}

Branched surfaces were first introduced by Williams in \cite{W}.  We refer the reader to \cite{FO}, \cite{O1} and \cite{O2} for the definitions of branched surface $B$ and associated $I$-bundle neighborhood $N(B)$.  We recall that a surface $S$, not necessarily compact, is {\sl carried} by a branched surface $B$ if $S$ is injectively immersed in $N(B)$ so that it is everywhere transverse to the $I$-fibers.  A minimal set of a foliation is carried by $B$ if each of its leaves is carried by $B$.  

It is standard to say that a foliation $\mathcal F$ is carried by $B$ if a Denjoy blow, $\mathcal F'$, of $\mathcal F$ along a single leaf $L$ results in a lamination $\Lambda=\mathcal F'\setminus j(L\times (0,1) $ such that all leaves of $\Lambda$ are carried by $B$.  (Note that if $B$ carries a foliation in this sense, then the regions complementary to $N(B)$ are necessarily products.) In this paper, however, we instead introduce and use the following definition.

\begin{definition}
A foliation $\mathcal F$ is {\sl carried by} a branched surface $B$ if 
\begin{enumerate}
\item the restriction of $\mathcal F$ to $N(B)$ is a foliation transverse to the $I$-fibers and tangent to $\partial_h N(B)$, and
\item the restriction of $\mathcal F$ to $M\setminus \mbox{int } N(B)$ is a product foliation that is transverse to the $I$-fibers of $N(B)$ along $\partial_v N(B)$ and tangent to $\partial_h N(B)$.
\end{enumerate}
\end{definition}
We refer the reader to \cite{GO} and \cite{Li} for a description of the splitting open of a branched surface and the corresponding ``splitting open'' of the $I$-bundle neighborhood $N(B)$.

\begin{thm}\label{closetofibration}
Let $\mathcal F$ be a taut, transversely oriented $C^0$ foliation of $M$ with only trivial holonomy.  Then $M$ fibers over $S^1$ and either 
\begin{enumerate}
\item $\mathcal F$ is measured, or else
\item $\mathcal F$ is a single Denjoy blow up of a minimal measured $C^0$ foliation.
\end{enumerate}
\end{thm}

\begin{proof} The proof proceeds by analyzing the three cases arising in the conclusion of Corollary~\ref{trivholo}.

Case~(1): If $\mathcal F$ is a fibering, of course $M$ is fibered.  It follows immediately that $\mathcal F$ satisfies conclusion~(1).

Case~(2): Since $\mathcal F$ does not contain an exceptional minimal set, Corollary~\ref{measured} can be applied.  From this it follows that $\mathcal F$ is transversely measured.  It remains to show that $M$ fibers over $S^1$.  Since $\mathcal F$ has a transverse measure, it is fully carried by a measured branched surface $B$, and since $\mathcal F$ is a foliation, the complementary regions of $B$ are products.

Case~(3): Then $\mathcal F$ is a single Denjoy blow up of a minimal foliation $\mathcal E$.  Denote the minimal set of $\mathcal F$ by $\Lambda$ and denote the smooth transverse flow used in the blow up by $\Phi$.  By Corollary~\ref{measured}, $\mathcal E$ is measured, thus (2) is satisfied.  It remains to show that $M$ fibers over $S^1$.

This measure on $\mathcal E$ determines a transverse measure on $\Lambda$ in the sense that any segment of a flow line of $\Phi$ can be given the measure it has when viewed as a transversal to $\mathcal E$.
Since $\Lambda$ has a transverse measure, it is fully carried by a measured branched surface $B$, and since the complementary regions of $\Lambda$ are $I$-bundles, the complementary regions of $B$ are products.  So Proposition~4.11 of \cite{O1} applies.

Hence, in each of the cases (2) and (3), $\mathcal F$ is fully carried by a transversely measured branched surface for which every complementary region is a product.  The space of positive measures on $B$ is an open cone in a vector space since it is the solution space of a system of homogenous linear equations.  Since there is a nontrivial real solution and the coefficients of these equations are integers, there is a nontrivial, positive, rational solution arbitrarily close to the real solution.  Any positive rational solution corresponds to an integral measure on $B$ and hence describes a surface $S$ (not necessarily connected) which is fully carried by $B$.  

Since any complementary region to $B$, and hence to $S$, is necessarily a product, it follows that $M$ cut open along $S$ is an $I$-bundle.  Letting $S_0$ be a (possibly, the) component of $S$, it follows that $M$ is a fiber bundle over $S^1$ with fiber $S_0$.
\end{proof}

Next we apply a theorem from \cite{KR4}.

\begin{thm}\label{measured near smooth} \cite{KR4} Suppose $\mathcal F$ is a transversely orientable $C^{1,0}$ measured foliation in $M$.  Then there is an isotopy of $M$ taking $\mathcal F$ to a $C^{\infty}$ measured foliation which is $C^0$ close to $\mathcal F$.  If $\Phi$ is a smooth flow transverse to $\mathcal F$, the isotopy may be taken to map each flow line of $\Phi$ to itself.
\end{thm}

\begin{thm}\label{measuredimpliesfibered} Fix $\epsilon>0$.  Let $\mathcal F$ be a taut, transversely oriented $C^{1,0}$ foliation of $M$ with only trivial holonomy.  Then $M$ fibers over $S^1$, and $\mathcal F$ is $O(\epsilon)$ $C^0$ close to a smooth fibering of $M$.  
\end{thm}

\begin{proof} 
By Theorem~\ref{closetofibration}, $M$ fibers over $S^1$, and $\mathcal F$ is a fibering, or a minimal measured foliation, or a Denjoy blow up of a minimal measured foliation.  If $\mathcal F$ is a fibering, there is nothing to prove.  And if $\mathcal F$ is measured and not a fibering, perform a single Denjoy blow up.

Thus, $\mathcal F$ is a single Denjoy blow up of a minimal measured $C^{1,0}$ foliation $\mathcal E$.  By Theorem~\ref{measured near smooth}, it suffices to show that $\mathcal F$ is $O(\epsilon)$ $C^0$ close to a $C^{1,0}$ fibering of $M$.  Denote the exceptional minimal set of $\mathcal F$ by $\Lambda$, and denote the smooth transverse flow used in the blow up by $\Phi$.  Let $B$ be a smoothly embedded, transversely oriented branched surface which is transverse to $\Phi$ and which fully carries $\mathcal F$.  In particular, $\Lambda$ lies in a regular $I$-fibered neighborhood $N$ of $B$, transverse to the $I$-fibers, with the $I$-fibers segments of the flow $\Phi$.  In addition, the regions complementary to $\text{int}(N)$ are sutured manifold products, and $\mathcal F$ restricts to a product foliation on these complementary regions.  Since $\Lambda$ has a transverse measure, so does $B$.

Now fix $\epsilon >0$.  Use $\Lambda$ to split $B$ open as much as is necessary to a smoothly embedded, measured branched surface $B_\text{split}$ so that the result of splitting open $N$ is a regular $I$-fibered neighborhood $N_\text{split}$ that has an $\epsilon$-flat $(\mathcal F, \Phi)$ flow box decomposition such that each flow box has horizontal boundary contained in $\partial_h N(B_\text{split})$.  We say $\mathcal G_0$ is a foliation of $N_\text{split}$ if in addition to the usual local product leaf structure, it is tangent to $\partial_h N_\text{split}$, transverse to $\partial_v N_\text{split}$, and everywhere transverse to $\Phi$.  By $\epsilon$-flatness, any foliation $\mathcal G_0$ of $N_\text{split}$ can be isotoped relative to $\partial N_\text{split}$ to be $\epsilon$ $C^0$ close to the restriction $\mathcal F_0$ of $\mathcal F$ to $N_\text{split}$.

Pick a rational measure $\mu$ fully carried by $B_\text{split}$, and let $\mathcal G_0$ be a $C^{1,0}$ fibering of $M$ determined up to isotopy by $(B,\mu)$.  Choose $\mathcal G_0$ to be everywhere transverse to $\Phi$ and so that its restriction to $N_\text{split}$ is a foliation of $N_\text{split}$.

Consider a component $\Sigma$ of the metric closure of the complement of $N_\text{split}$.  Both $\mathcal F$ and $\mathcal G_0$ restrict to product foliations on $\Sigma$, and hence to foliations by circles on each component $A_i$ of $\partial_v N_\text{split}$.  Again appealing to the $\epsilon$-flatness of $N(B_\text{split})$, there is an $O(\epsilon)$ $C^0$ small isotopy of $M$ in a neighborhood of $\cup_i A_i$ taking $\mathcal F$ to $\mathcal F_0$ so that if $F$ and $G$ are leaves of the restrictions of $\mathcal F_0$ and $\mathcal G_0$ respectively to $\Sigma$, then either $F= G$ in a neighborhood of $\partial F=\partial G$ or $\partial F\cap\partial G=\emptyset$.  Now let $\mathcal G$ be the foliation obtained by letting $\mathcal G$ coincide with $\mathcal G_0$ on $N_\text{split}$ and coincide with $\mathcal F_0$ on each component $\Sigma$.  By construction, $\mathcal G$ is a $C^{1,0}$ fibering of $M$.  Moreover, $\mathcal G$ is $O(\epsilon)$ close to $\mathcal F$ since $\mathcal G_0$ is $O(\epsilon)$ close to $\mathcal F$ on $N_\text{split}$ and coincides with $\mathcal F_0$ on each $\Sigma$.
\end{proof}

Recall the statement of Tischler's Theorem.

\begin{thm} \cite{tischler} \label{tischler}A transversely oriented, $C^\infty$ measured foliation $\mathcal F$ of $M$ can be $C^\infty$ approximated by a smooth fibering $\mathcal G $ of $M$ over $S^1$.
\end{thm}

As a corollary to Theorem~\ref{measuredimpliesfibered}, we thus have Tischler's Theorem for 3-manifolds $M$.

Finally, we highlight the special case that all leaves in some minimal set of $\mathcal F$ are planes.  

\begin{thm}\cite{GICM}\label{gicm0} If $\lambda$ is an essential lamination of the closed 3-manifold $M$ by planes, then $M$ is the 3-torus.
\end{thm}

As an immediate corollary of Theorem~\ref{gicm0} and its proof in \cite{GICM}, we conclude the following.

\begin{cor} \label{gicm} Let $\mathcal F$ be a Reebless $C^0$ foliation of a closed 3-manifold $M$.  Suppose $\mathcal F$ contains a minimal set all of whose leaves are planes.  Then all leaves of $\mathcal F$ are planes and $M=T^3$.  
\end{cor}

\begin{proof} Gabai's theorem in \cite{GICM} is stated for essential laminations.  We translate as follows.  Let $\Lambda$ be a minimal set of $\mathcal F$ all of whose leaves are planes.  Since $\mathcal F$ is Reebless, $\Lambda$ is essential.  In fact, if $B$ is any essential branched surface which fully carries $\Lambda$, then any component of $\partial_h N(B)$ is necessarily a disk, and, since $M$ is irreducible, the complementary regions of $B$ are necessarily homeomorphic to product regions $D^2\times I$.  Since leaves of $\mathcal F$ are $\pi_1$-injective, all leaves of $\mathcal F$ must be planes.  As described in \cite{GICM}, it follows that $M=T^3$.
\end{proof}

The $C^0$ foliation theory that is used for the main result in this paper is the following corollary.

\begin{cor}\label{allplanes}
Let $\mathcal F$ be a taut, transversely oriented $C^{1,0}$ foliation of a closed 3-manifold $M$.  Suppose $\mathcal F$ contains a minimal set all of whose leaves are planes.  Then $M=T^3$ and $\mathcal F$ is $C^0$ close to a smooth fibering of $M$.  
\end{cor}

\begin{proof} By Theorem~\ref{gicm}, all leaves of $\mathcal F$ are planes and $M=T^3$.  Necessarily, therefore, $\mathcal F$ has only trivial holonomy, and so by Theorem~\ref{measuredimpliesfibered}, $\mathcal F$ is $C^0$ close to a smooth fibering of $M$.
\end{proof}

\section{Propagation} \label{propagate}

In this section, we describe how to extend the smooth confoliation on $V$ to a smooth contact structure on $M$ which is $O(\epsilon)$ close to $\mathcal F$, and hence prove our main result, Theorem~\ref{maincor}.

Begin by recalling the propagation technique introduced in \cite{KR2}.  The starting point is a decomposition of a manifold into two codimension-0 pieces.  Roughly speaking, one piece, $V$, has a contact structure, while the other piece, $W$, has a foliation.  As long as the contact structure dominates the foliation along $\partial_vV=\partial_vW$ and every point of $W$ can be connected to $V$ by a path in a leaf, the contact structure can be propagated throughout $W$.  
 
For completeness we include the formal definitions of these concepts and the main theorem of \cite{KR2} stated in the $C^0$ setting.

\begin{definition}[Definition~6.4 of \cite{KR2}] \label{compatible1} Let $M$ be a closed oriented 3-manifold with smooth flow $\Phi$.  Suppose that $M$ can be expressed as a union 
$$M=V\cup W,$$ where $V$ and $W$ are smooth 3-manifolds, possibly with corners, such that $\partial V=\partial W$.  We say that this decomposition is {\sl compatible with the flow $\Phi$} if $\partial V$ (and hence $\partial W$) decomposes as a union of compact subsurfaces $\partial_v V\cup \partial_h V$, where $\partial_v V$ is a union of flow segments of $\Phi$ and, $\partial_h V$ is transverse to $\Phi$.  Let $U$ be a preferred regular neighborhood of the union of the horizontal 2-cells and the vertical 1-cells of $\partial V$.  
Suppose this decomposition is compatible with $\Phi$, $V$ admits a smooth confoliation $\xi_V$, and $W$ admits a $C^{1,0}$ foliation $\mathcal F_W$.  Suppose that $U$ is smoothly foliated with a foliation $\mathcal F_U$ which smoothly agrees with $\mathcal F_W$ where they meet.  We say that $(V,\xi_V)$ is {\sl $\Phi$-compatible} with $(W,\mathcal F_W)$, and that $M$ admits a positive $(\xi_V,\mathcal F_W,\Phi)$ decomposition, if the following are satisfied:
\begin{enumerate}
\item $\xi_V$ and $\mathcal F_W$ are (positively) transverse to $\Phi$ on their domains of definition,

\item each of $\mathcal F_W$ and $\xi_V$ is tangent to $\partial_h V$,

\item $\xi_V = T\mathcal F_U$ on $\overline U \cap V$, 

\item $\xi_V$ is a contact structure on $V\backslash U$, and 

\item $\chi_{\xi_V}<\chi_{T\mathcal F_W}$ on $(\partial_v V) \backslash\overline U$, when viewed from outside $W$.

\end{enumerate}
A foliation $\mathcal F_W$ is {\sl $V$-transitive} if every point in $W$ can be connected by a path in a leaf of $\mathcal F$ to a point of $V$.
\end{definition}

Let $V_{\gamma_i}({\tau_i})$ be the spanning collection of attracting neighborhoods constructed in Theorem~\ref{metriccreateholonbhds}.  Then $M$ can be decomposed by setting $V= V_{\gamma_1}(\tau_1) \cup...\cup V_{\gamma_n}(\tau_n)$ and letting $W$ be the closure of the complement of $V$.  Since $V$ is the union of a spanning collection of neighborhoods, $\mathcal F'_W$, the restriction of $\mathcal F'$ to $W$, is $V$-transitive.

The confoliation $\xi_V$ constructed in Corollary~\ref{extend} satisfies Conditions (1)--(5) , and we can therefore apply Theorem~6.10 of \cite{KR2}.

\medskip

\begin{thm}[Theorem~6.10 of \cite{KR2}] \label{main1} If $M$ admits a positive $(\xi_V,\mathcal F_W,\Phi)$ decomposition such that $\mathcal F_W$ is {\sl $V$-transitive} and $\xi_V$ is $\epsilon$ $C^0$ close to $\mathcal F_W \cap \partial_vW$, then $M$ admits a smooth positive contact structure $\xi^+$ which agrees with $\xi_V$ on $V$ and is $\epsilon$ $C^0$ close to $\mathcal F_W$ on $W$.  The analogous result holds if $M$ admits a negative $(\xi_{V'},\mathcal F_{W'},\Phi)$ decomposition, yielding a smooth negative contact structure $\xi^-$.  If $M$ admits both a positive $(\xi_V,\mathcal F_W,\Phi)$ decomposition and a negative $(\xi_{V'},\mathcal F_{W'},\Phi)$ decomposition, then these contact structures $(M,\xi^+)$ and $(-M,\xi^-)$ are weakly symplectically fillable and universally tight.
\end{thm}

Our main theorem can now be proved.
\medskip

\noindent{\bf Theorem~\ref{maincor}.} {\it Any taut transversely oriented $C^{1,0}$ foliation on a closed oriented 3-manifold $M\ne S^1\times S^2$ can be $C^0$ approximated by both a positive $\xi^+$ and a negative $\xi^-$ smooth contact structure.  These contact structures $(M,\xi^+)$ and $(-M,\xi^-)$ are weakly symplectically fillable and universally tight.}

\medskip

\begin{proof}[Proof of Theorem~\ref{maincor}]
 Given the earlier results in this paper on approximating taut foliations by smoother taut foliations, the methods developed for introducing holonomy in minimal sets, and the construction of approximating contact structures in attracting holonomy neighborhoods, Theorem~\ref{maincor} follows directly from Theorem~\ref{main1}.  Since this is a long line of implications and constructions, we assemble and summarize the steps now.

The first step is to show that it suffices to restrict attention to the case that $\mathcal F$ is a taut, transversely oriented $C^{\infty,0}$ foliation which is not a fibering and whose every minimal set contains a leaf which is not homeomorphic to $\mathbb R^2$.

Consider first the case that $\mathcal F$ is a taut, transversely oriented $C^{1,0}$ foliation on $M$.  By Theorem~\ref{calegarikr}, $\mathcal F$ can be $C^0$ approximated by a taut, transversely oriented $C^{\infty,0}$ foliation.  Next if $\mathcal F$ contains a minimal set all of whose leaves are planes, then by Corollary~\ref{allplanes}, it can be $C^0$ approximated by a smooth fibering.  Finally, if $\mathcal F$ is a $C^{\infty,0}$ fibering, then by Theorem~\ref{blowupthm}, it can be $C^0$ approximated by a taut, transversely oriented $C^{\infty,0}$ foliation which is obtained by Denjoy blow up.  

If $\mathcal F$ is minimal, set $\Lambda_1=\mathcal F$.  Otherwise, let $ \Lambda_1,...,\Lambda_r$ denote the exceptional minimal sets of $\mathcal F$, and let $[L_1],...,[L_s]$ denote the isotopy classes of compact leaves of $\mathcal F$.  Apply Corollary~\ref{finitespan} to obtain a spanning collection of pairwise disjoint holonomy neighborhoods $V'_{\gamma_1}(\tau_1,A_1),...,V'_{\gamma_{r+s}}(\tau_{r+s},A_{r+s}).$ Let $V'$ denote their union.

For each $i, 1\le i\le n$, let $R'_i=R_{\gamma_i}(\sigma_i,A_i)$, and set $R'=\cup_i R_i$.  For each $i, 1\le i\le n$, fix a smooth open neighborhood $N_{R_i'}$ of $R_i'$ in $V_i'$.  Choose each $N_{R_i'}$ small enough so that its closure, $\overline{N_{R_i'}}$, is a closed regular neighborhood of $R_i'$.  Let $N_{R'}$ denote the union of the $N_{R_i'}$.

By Lemma~\ref{stronglycompatibleF}, $\mathcal F$ can be $C^0$ approximated by a taut, transversely oriented $C^{\infty,0}$ foliation which is strongly $(V',P)$ compatible for some choice of product neighborhood $(P,\mathcal P)$ of $(V';N_{R'})$.  Hence, it suffices to restrict attention to the case that $\mathcal F$ is a taut, transversely oriented, strongly $(V',P)$ compatible, $C^{\infty,0}$ foliation which is not a fibering and whose every minimal set contains a leaf which is not homeomorphic to $\mathbb R^2$.  We now do so.  In particular, $\mathcal F=\mathcal P$ on $N_{R'}$ and $\mathcal F$ is $x$-invariant in the $(x,y,z)$ coordinates given by $P$.

Put the product metric on each component $P_i=[-1,1]\times S^1\times [-1,1]$ of $P$, as described in Definition~\ref{metricpn}, and let $g_0$ denote the resulting metric on $P$.  Fix a Riemannian metric $g=g(P)$ on $M$ which restricts to $g_0$ on $P$.  Fix $\epsilon>0$.

Now apply Theorem~\ref{metriccreateholonbhds} to obtain a taut, transversely oriented $C^{\infty,0}$ foliation $\mathcal G$ that is $\epsilon$ $C^0$ close to $\mathcal F$, $V'$-compatible with $\mathcal F$, and strongly $(V',P)$ compatible, and a finite set of pairwise disjoint attracting neighborhoods $$V_{\gamma_1}(\tau_1,A_1),...,V_{\gamma_m}(\tau_m,A_m), $$ $m\ge n$, for $\mathcal G$ such that the $V_{\gamma_i}$ are $\epsilon$-flat, $\epsilon$-horizontal, and if $V$ denotes their union, 
$\mathcal G$ is $V$-transitive.

By Theorem~\ref{C0V}, there are a regular neighborhood $N_v \subset V$ of the vertical edges of $\partial Q$ in $V$ and smooth approximating confoliations $\xi_V^{\pm}$ defined on $V$ where the main properties are 
\begin{enumerate}

\item $\xi_V^{\pm}=T\mathcal G$ on $\overline{N_h\cup N_v}$ and is contact at all other points of $V$,

\item $\xi_V^+$ dominates $\mathcal G$ along $\partial_v V$, with the domination strict outside $\overline{N_h\cup N_v}$, and

\item $\xi_V^-$ is dominated by $\mathcal G$ along $\partial_v V$, with the domination strict outside $\overline{N_h\cup N_v}$.
\end{enumerate}

Let $W$ denote the closure of the complement of $V$ and let $\mathcal G_W$ denote the restriction of $\mathcal G$ to $W$.  Applying Theorem~\ref{main1} to the $(\xi_V^+,\mathcal G_W, \Phi)$ decomposition of $M$ yields a positive contact structure $\xi^+$ on $M$ which is $O(\epsilon)$ close to $\mathcal G$.  By symmetry, there is a negative contact structure $\xi^-$ on $M$ which is $O(\epsilon)$ close to $\mathcal G$.  These contact structures are weakly symplectically fillable and universally tight.  Since $\mathcal G$ is $\epsilon$ close to $\mathcal F$, each of $\xi^{\pm}$ is $O(\epsilon)$ close to $\mathcal F$.
\end{proof}


\begin{thebibliography}{}

\bibitem{bowden}
J.\ Bowden, \textit{Approximating $C^0$-foliations by contact structures}, (preprint).

\bibitem{goodman}
C.\ Camacho and A.\ Neto, (translated by S.\ E.\ Goodman), \textit{Geometric theory of foliations}, Birkh\"auser, Boston, 1985.

\bibitem{calegari} D.\ Calegari, \textit{Leafwise smoothing laminations}, Algebr.  Geom.  Topol.  \textbf{1} (2001), 579--585.

\bibitem{CC} A.\ Candel and L.\ Conlon, \textit{Foliations I}, A.M.S.  Graduate Studies in Mathematics \textbf{23}, 2000.

\bibitem{DL} O.\ Dasbach and T.\ Li, \emph{Property P for knots admitting certain Gabai disks}, Topology and its Applications \textbf{142} (2004), 113--129.

\bibitem{Di} P.\ Dippolito, \textit{Codimension one foliations of closed manifolds}, Ann.\ of Math, {\bf (2)}, 1978, pp.\ 403--453.

\bibitem{DR} C.  \ Delman and R.\ Roberts, \emph{Alternating knots satisfy Strong Property P}.  
Comment.  Math.  Helv.  \textbf{74} (1999), no.  3, 376--397.  

\bibitem{ET}
Y.\ Eliashberg and W.\ Thurston, \textit{Confoliations}, University Lecture Series \textbf{13}, Amer.\ Math.\ Soc., Providence, 1998.  

\bibitem{FO}
W.\ Floyd and U.\ Oertel, \textit{Incompressible branched surfaces via branched surfaces}, Topology \textbf{23} (1984), 117--125.



\bibitem{G1} D.\ Gabai, \emph{Foliations and the topology of 3-manifolds}, J.  Differential Geometry.  \textbf{18} (1983), 445-503.

\bibitem{g1} D.\ Gabai, \emph{Foliations and Genera of Links}, Topology \textbf{23} (1984), 381--394.

\bibitem{g2} D.\ Gabai, \emph{Genera of the Alternating Links}, Duke Math.  J.  \textbf{53} (1986), 677--681.

\bibitem{g3} D.\ Gabai, \emph{Detecting Fibred Links in S3}, Commentarii Math.  Helvetici \textbf{61} (1986), 519--555.

\bibitem{G2} D.\ Gabai, \emph{Foliations and the topology of 3-manifolds II} J.  Differential Geometry.  \textbf{26} (1987), 461--478.

\bibitem{G3} D.\ Gabai, \emph{Foliations and the topology of 3-manifolds III.} J.  Differential Geometry.  \textbf{26} (1987), 479-536.

\bibitem{GICM} D.\ Gabai, \emph{Foliations and 3-manifolds}, Proc.  I.C.M.  Vol.  I, II (Kyoto, 1990), 609--619, Math.  Soc.  Japan, Tokyo, 1991.

\bibitem{Ga1}
D.\ Gabai, \textit{Taut foliations of 3-manifolds and suspensions of $S^1$}, Annales de l'institut Fourier, tome {\bf 42}, no.  1--2 (1992), pp.  193--208.

\bibitem{Ga93}
D.\ Gabai, \textit{Problems in foliations and laminations}, in Geometric Topology (Athens, GA, 1993) AMS/IP Stud.  Adv.  math.  2.2, pp.  1?-33.

\bibitem{Ga2}
D.\ Gabai, \textit{Combinatorial volume preserving flows and taut foliations}, Comment.  Math.  Helv., tome {\bf 75} (2000), pp.  109--124.

\bibitem{GO}D.\ Gabai and U.\ Oertel, {\it Essential Laminations in 3-Manifolds}, Ann.  Math.  {\bf 130}, (1989), 41--73.

\bibitem{GP}
V.\ Guillemin and A.\ Pollack, \textit{Differential Topology}, Prentice-Hall 1974.

\bibitem{Hass} J.\ Hass, \textit{Minimal surfaces in foliated manifolds}, Comment.  Math.  Helv.  \textbf{61} (1986), 1--32.

\bibitem{HH}
G.\ Hector and U.\ Hirsch, \textit{Introduction to the geometry of foliations, Part B}, Aspects of Mathematics, Vol.  E3 (1983).

\bibitem{Im}
H.\ Imanishi, \textit{On the theorem of Denjoy-Sacksteder for codimension one foliations without holonomy}, J.  Math.  Kyoto University \textbf{14-3} (1974), pp.  607--634.

\bibitem{KRo} T.\ Kalelkar and R.\ Roberts, \textit{Taut foliations in surface bundles with multiple boundary components}, Pacific J.  Math.  {\bf 273} (2015), pp.  257--275.

\bibitem{KR2}
W.\ H.\ Kazez and R.\ Roberts, \textit{Approximating $C^{1,0}$-foliations}, Geometry \& Topology Monographs (2015) (to appear).

\bibitem{KR4}
W.\ H.\ Kazez and R.\ Roberts, \textit{Flow box decompositions}, (preprint).  



\bibitem{KR5}
W.\ H.\ Kazez and R.\ Roberts, \textit{Taut foliations}, (preprint).

\bibitem{Lee} J.\ M.\ Lee, \textit{Introduction to smooth manifolds}, Graduate Texts in Mathematics \textbf{218} (2003).

\bibitem{Li}
T.\ Li, {\it Laminar branched surfaces in 3-manifolds}, Geom.\ Topol.\ {\bf 6} (2002), pp.  153--194 (electronic).

\bibitem{Li2}
T.\ Li, \text{Commutator subgroups and foliations without holonomy}, Proc.  Amer.  Math.  Soc.  {\bf 130} (2002), no.  8, pp.  2471--2477.

\bibitem{Li3} T.\ Li, \emph{Boundary train tracks of laminar branched surfaces}.  Proc.  of Symposia in Pure Math., \textbf{71}, AMS (2003), 269--285.

\bibitem{LR} T.\ Li and R.\ Roberts, {\it Taut foliations in knot complements}, Pacific J.  Math.  {\bf 269} (2014), pp.  149--168.

\bibitem{Mo}
E.\ E.\ Moise, \textit{Affine structures in 3-manifolds.  V.  The triangulation theorem and Hauptvermutung}, Annals of Math.  Second Series {\bf 56} (1952), pp.  96--114.

\bibitem{munkres}
J.\ R.\ Munkres, \textit{Elementary differential topology}, Annals of Math.  Studies \textbf{54}, 1966.

\bibitem{OS}
P.\ Ozsv\'ath and Z.\ Szab\'o, {\it Holomorphic disks and genus bounds}, Geometry and Topology \textbf{8} (2004), 311--334.

\bibitem{Novikov}
S.  \ P.\ Novikov, \emph{Topology of Foliations}, Trans.  Moscow Math.  Society (1965), 268--304.

\bibitem{O1} U.\ Oertel, \textit{Incompressible branched surfaces}, Invent.  Math.  \textbf{76} (1984), 385--410.

\bibitem{O2} U.\ Oertel, \textit{Measured laminations in 3-manifolds}, Trans.  A.M.S.  \textbf{305(2)} (1988), 531--573.


\bibitem{rado}
T.\ Rad\'o, \textit{\"Uber den begriff der Riemannschen fl\"ache}, Acta Univ.  Szeged, \textbf{2} (1924-1926), 101--121.

\bibitem{R} R.\ Roberts, \emph{Constructing taut foliations}, Comment.  Math.  Helv.  \textbf{70} (1995), 516--545.

\bibitem{R1} R.\ Roberts, \emph{Taut foliations in punctured surface bundles, I}.  Proc.  London Math.  Soc.  (3) \textbf{82} (2001), no.  3, 747--768.

\bibitem{R2} R.\ Roberts, \emph{Taut foliations in punctured surface bundles, II}.  Proc.  London Math.  Soc.  (3) \textbf{83} (2001), no.  2, 443--471.

\bibitem{S} R.\ Sacksteder, \textit{Foliations and pseudogroups}, Amer.  J.  Math.  \textbf{87} (1965), 79--102.

\bibitem{SS} R.\ Sacksteder and A.  Schwartz, \textit{Limit sets for foliations}, Ann.  Inst.  Fourier \textbf{15} (1965), 201--214.

\bibitem{Siebenmann}
L.\ C.\ Siebenmann, \textit{Deformations of homeomorphisms on stratified sets}, Comm.  Math.  Helv.  \textbf{47} (1972), 123--163.

\bibitem{So}
V.\ V.\ Solodov, \textit{Components of topological foliations}, Mat.  Sb.  \textbf{119(161)} (1982), 340--354.
English translation in Math.  USSR Sb.  \textbf{47} (1984), 329--343.

\bibitem{Sullivan} D.\ Sullivan, \textit{Cycles for the dynamical study of foliated manifolds and complex manifolds}, Inventiones.  Math.  \textbf{36} (1976), 225--255.  

\bibitem{thurston} 
W.\ P.\ Thurston, \textit{A Norm for the homology of 3-manifolds}, Memoirs of the AMS \textbf{59(339)}, 99--130.

\bibitem{tischler}
D.\ Tischler, \textit{On fibering certain foliated manifolds over $S^1$}, Topology \textbf{9} (1970), 153--154.

\bibitem{vogel}
T.\ Vogel, \textit{Uniqueness of the contact structure approximating a foliation}, arXiv:1302.5672v2


\bibitem{W} R.\ Williams, \emph{Expanding attractors}, Inst.  Hautes \'Etudes Sci.  Publ.  Math.  \textbf{43} (1974), 169--203.

\end{thebibliography}
\end{document}